# MODIFIED CONVEX HULL PRICING FOR POWER MARKETS WITH PRICE-SENSITIVE LOAD

Vadim Borokhov[1]

LLC "En+development",
ul.Vasilisi Kozhinoi 1, Moscow, 121096, Russia

**Abstract**

We consider a general power market with price-sensitive consumer bids and non-convexities originating from supply (start-up and no-load costs, nonzero minimum output limits of generating units, etc.) and demand. The convex hull (minimum-uplift) pricing method produces the set of power prices that minimizes the total uplift payments to the market players needed to compensate their potential profits lost by accepting the centralized dispatch solution. All opportunities to supply (consume) any other output (consumption) volumes allowed by market player individual operational constraints are considered as foregone in the convex hull pricing method. We modify the convex hull pricing algorithm by defining for each market player a modified individual feasible set that is utilized in the lost profit calculation. These sets are based on the output (consumption) volumes that are economically and technologically feasible in the centralized dispatch. The new pricing method results in the generally different set of market prices and lower (or equal) total uplift payment compared to the convex hull pricing algorithm.

## 1. Introduction

Deregulation of the electric power markets led to a free pricing for power realized through the bilateral trade and/or the centrally coordinated price auction. The latter usually has the form of a bid-based security-constrained optimization problem, which produces both the generating unit dispatch and load schedule as well as the marginal price for power [1]-[3]. If the corresponding optimization problem is convex, then the marginal price is an equilibrium price: if the market price is set equal to the marginal price, then no market player, acting as a price-taker, has incentives to change its dispatch/load schedule within the player's individual (private) feasible set specified in the problem. (If there is a set of marginal prices, then there exists a marginal price that is also an equilibrium price.) However, if some of the power output cost components are not reflected in the generator bids, application of the marginal pricing doesn't guarantee non-confiscatory pricing. This may occur, for example, in case of the convex economic dispatch problem, which has fixed unit commitment statuses and, therefore, doesn't account for non-convexities in the generator output cost functions such as nonzero minimum output limits and fixed cost components of power production. If the optimization problem is not convex, then an equilibrium price may not exist and some other pricing mechanism is needed to ensure economic stability of the centralized dispatch outcome and non-confiscatory market pricing for power [4]. The non-convexities may originate from the flexible demand side due to discrete and minimum power consumption levels [5] as well as the supply side because of the fixed and start-up costs, minimum output limits, integral commitment decision variables, minimum up/down times, etc. It is known that the absence of an equilibrium market price for power in uninode single-period markets with non-convexities, originating from fixed cost component of power output and/or nonzero minimum output limits, can be traced to nonintersection of the aggregate demand and the aggregate supply curves due to gaps (discontinuities) in the aggregate supply curve.

---

[1] E-mail: vadimbor@yahoo.com. The views expressed in this paper are solely those of the author and not necessarily those of LLC "En+development".

A number of pricing methods for the centrally coordinated markets with non-convexities has been developed, [6]-[22]. In [6] the new services (units being online) and the corresponding prices are introduced in addition to electric power by setting the integral status variables to their optimal values obtained from the centralized dispatch solution. However, the price for being online can be negative and the resulting pricing produces zero profit of generators and is close to pay-as-bid pricing. If the negative prices are discarded in order to allow generators to retain their profits, then the resulting pricing does not achieve the competitive equilibrium. In [7], [8] this scheme was modifies to generate more stable prices by adding extra constraints to the reformulated optimization problem that fix certain continuous variables at their optimal values as well.

The nonlinear (discriminatory) pricing for power with market player specific prices was proposed in [9]-[12]. The nonlinear pricing in the form of the generalized uplift functions that includes generators as well as consumers in the lost profit compensation and ensures zero net uplift at the market was introduced in [13], [14], [15]. Authors in [16] proposed the minimum zero-sum uplift pricing approach that increases the price above marginal cost and transfers all the additional payments (that the profitable suppliers receive as a result of the price increase) to the unprofitable suppliers to make them whole in the form of internal zero-sum uplifts. In [17] a primal-dual approach was proposed to find the market prices that minimize the social welfare reduction due to schedules inconsistency and ensure non-negative generator profits. Such an approach, however, doesn't guarantee the opportunity cost compensation for generators; therefore, competitive equilibrium at the centralized dispatch solution is not achieved. In [18] a semi-Lagrangian relaxation approach was developed to compute a uniform market price that produces the same solution as the original centralized dispatch problem while ensuring that no supplier incurs losses. A zero-sum uplift pricing scheme that aims to minimize the maximum contribution to the financing of the uplifts in a model where both suppliers and buyers place bids was suggested in [19]. If the demand is fixed, such an approach results in the market price equal the maximum average cost of the producers.

The convex hull pricing method proposed in [20]-[22] is a pricing mechanism with uniform market price that minimizes the corresponding total uplift payment needed to ensure economic stability of the centralized dispatch outcome. Each market player is assumed to have an opportunity to supply (consume) any power volume allowed by its individual operational constraints and is compensated for an opportunity to supply (consume) the individually feasible volume with the profit higher than that obtained when following the centralized dispatch at a given value of the market price. In this framework, the total lost profit to be compensated through the uplift payments coincides with the duality gap that emerges after the Lagrangian relaxation procedure is applied to the power balance constraint [20]-[22]. Thus, the convex hull pricing implies that given the market price each market player receives profit that is equal to the maximum value of its profit function on the set specified by the market player individual operational constraints only. Such payments ensure the economic stability of the centralized dispatch outcome. In this approach, the market price is the only exogenous parameter that encodes information on the external power system and is utilized to calculate the lost profit of the market player. If the centralized dispatch optimization problem is convex, then the convex hull pricing mechanism results in the standard marginal prices. In [5] the convex hull pricing and generalized uplift approaches were applied to the market with both demand-side and supply-side non-convexities.

The total uplift payment includes transmission congestion and loss revenues that determine the revenue adequacy of financial transmission rights [23]. The convex model that solves unit commitment problem for practical instances and can be readily employed for the convex hull pricing method was developed in [24], [25]. Detailed analysis of various aspects of the convex hull pricing method was presented in [26].

Clearly, the uplift payments distort the uniform market pricing and decrease transparency of the market. Therefore, it is all-important to reduce these payments. In [27] it was proposed to reduce the total uplift payment at the expense of having the new redundant constraint introduced in the optimization problem. The new constraint depends linearly on the unit status variables of all generators. Therefore, to obtain decoupling of the relaxed problem into a set of independent decentralized dispatch problems the relaxation procedure needs to be applied to the newly introduced constraint as well as the power balance constraint. This leads to the introduction of a new service (a unit being online) and the associated price in addition to the market price for power.

In [28] the market opportunities of generators foregone by accepting the centralized dispatch solution were studied for the case of single-period uninode power system with fixed demand and zero minimum output limits of all the generating units. For each generator, the modified individual feasible set corresponding to a set of technologically and economically feasible output volumes was formulated using the additional (redundant) constraints depending on the output volumes of that generator only. Hence, no new products/services and the associated prices are introduced in this approach. The algorithm excludes the output volumes that are technologically and/or economically infeasible from the uplift payment calculation: given the market price each market player receives profit that is equal to the maximum value of its profit function on the modified individual feasible set, which corresponds to the set of technologically and economically feasible output volumes of the market player. The resulting modified convex hull pricing approach yields the total uplift payment less or equal than that obtained in the convex hull pricing method. If the centralized dispatch optimization problem is convex, then application of the new method produces the standard marginal prices. Thus, contrary to the convex hull pricing method, the modified approach implies that, in addition to the market price, ISO informs each market player about the new constraints (which generally differ for the different market players) that together with the market player individual operational constraints specify the modified individual feasible set of the market player. Such a pricing approach can be viewed as augmenting the producer revenue function (before the uplift) to include the nonlinear terms representing the penalty functions for the new constraints without modifying the individual feasible set of the market player. From this viewpoint, the approach proposed in [28] is a combination of nonlinear pricing method (with nonlinear terms vanishing at the solution) and minimum-uplift linear pricing principle.

It was shown in [28] that the convex hull pricing method and the modified approach may produce the different outcomes in uninode single-period models with fixed load and zero minimum output limits of all units only if there is at least one generator with decreasing average output cost function up to the volumes exceeding demand (which also implies that the generator's maximum output limit exceeds demand). However, the consideration in [28] was limited to single-period uninode power systems with fixed demand and zero minimum output limits of the generators, i.e. the only source of non-convexity was the fixed cost of power output.

In this paper, we extend the modified convex hull pricing method [28] to multi-period power systems with price-sensitive load bids and both demand-side and supply-side non-convexities. We also allow for nonzero minimum output limits of generating units. For simplicity, we restrict our consideration to the uninode market model describing a power system without transit power losses and network constraints. In our study, all the market players are assumed to act as price-takers (in economic sense). We also outline extension of the method to multi-node systems and present an example illustrating application of the proposal for two-node case with the network constraint.

The paper is organized as follows. In Section 2 we briefly review the convex hull pricing method for both fixed and price-sensitive load cases. In Section 3 we define the modified individual feasible sets for the market players corresponding to opportunities available for a market player in the centralized power market. The proposed modified convex hull pricing method is introduced for both fixed and price-sensitive load cases in Sections 4 and 5, respectively. In Section 6 we study properties of the market prices and the uplift payments resulting from the new method, we also consider examples illustrating application of the new approach (including the cases of the demand-side non-convexity, multi-node power system, and time-coupled multi-period market) and make comparisons with implications of the convex hull pricing method. The main results are summarized in Section 7.

## 2. Convex hull pricing

Consider the centrally dispatched $T$-period uninode power market with fixed demand vector $\mathbf{d} = (d^1,...,d^T)$, $\mathbf{d} \in R_{\geq 0}^T$, and $n$, $n < +\infty$, generating units. For a unit $i$, $i \in I$, $I = \{1,...,n\}$, let $u_i^t$ be the binary commitment status variable taking values in the set $\{0,1\}$ (with 0 for the unit with status "OFF" and 1 for the unit with status "ON") and $g_i^t$ be the unit output volume at time interval $t$, $t \in T$, where $T$ is a set of time periods of a given market planning horizon. Denote as $\mathbf{u}_i$, $\mathbf{u}_i \in \{0,1\}^T$, the commitment vector and as $\mathbf{g}_i$, $\mathbf{g}_i \in R_{\geq 0}^T$, the power output vector of the unit $i$. Let $X_i$ be the individual feasible set of a generator $i$, $i \in I$. We note that in the general case the set $X_i$ depends on the initial conditions, which may include statuses and outputs volumes in the present market planning cycle as well as the commitment and power output volume in the previous market planning cycles. For simplicity, we do not explicitly indicate the initial conditions as they can be easily incorporated in the consideration. The set $X_i$ is assumed to be compact, $\forall i \in I$. Let a unit $i$ bids the cost function $C_i(\mathbf{x}_i)$, where $\mathbf{x}_i = (\mathbf{u}_i, \mathbf{g}_i)$. For any fixed $\mathbf{u}_i$, $C_i(\mathbf{x}_i)$ viewed as a function of $\mathbf{g}_i$ is assumed to be the non-negative continuous convex function in its corresponding domain (which may depend on the value of $\mathbf{u}_i$). The centralized dispatch optimization problem (which we refer to as the primal problem) with optimization (decision) variables $\mathbf{x} = (\mathbf{x}_1,...,\mathbf{x}_n)$ has the form

$$v = \min_{\substack{\mathbf{x}, \\ \mathbf{x}_i \in X_i, \forall i \in I \\ \sum_{i \in I} \mathbf{g}_i = \mathbf{d}}} \sum_{i \in I} C_i(\mathbf{x}_i), \quad (1)$$

where $v$ is the total cost to meet demand. Let us denote by $\Omega$ the feasible set specified by constraints in (1), the set $\Omega$ is assumed to be nonempty and compact. Together with continuity of $C_i(\mathbf{x}_i)$ in $\mathbf{g}_i$, this ensures existence of minimum in (1). Let $\mathbf{x}^* = (\mathbf{x}_1^*,..,\mathbf{x}_n^*)$ with $\mathbf{x}_i^* = (\mathbf{u}_i^*, \mathbf{g}_i^*)$ denote an optimal point of (1). If (1) has multiple optimal points, then $\mathbf{x}^*$ denotes any of them.

The decentralized dispatch problem for a generating unit $i$ at a given vector of the market prices $\mathbf{p} = (p^1,..,p^T)$ has the form

$$\pi_i^{prod}(\mathbf{p}) = \max_{\mathbf{x}_i \in X_i} \pi_i^{prod}(\mathbf{p}, \mathbf{x_i}), \text{ with } \pi_i^{prod}(\mathbf{p}, \mathbf{x_i}) = \mathbf{p}^T \mathbf{g}_i - C_i(\mathbf{x}_i). \quad (2)$$

We note that $\pi_i^{prod}(\mathbf{p})$ is a proper convex function with the domain $dom[\pi_i^{prod}(\mathbf{p})] = R^T$. Now we are ready to formulate the dual of the primal optimization problem (1):

$$v^D = \max_{\mathbf{p} \in R^T} \min_{\mathbf{x}, \mathbf{x}_i \in X_i, \forall i \in I} \left[ \mathbf{p}^T \left( \mathbf{d} - \sum_{i \in I} \mathbf{g}_i \right) + \sum_{i \in I} C_i(\mathbf{x}_i) \right] = \max_{\mathbf{p} \in R^T} \left( \mathbf{p}^T \mathbf{d} - \sum_{i \in I} \pi_i^{prod}(\mathbf{p}) \right). \quad (3)$$

As (3) is the unconstrained maximization problem of the concave function, its set of maximizers, which we denote as $P^+$, is given by the solutions to $\mathbf{d} \in \sum_{i \in I} \partial \pi_i^{prod}(\mathbf{p})$. It is straightforward to see that $P^+$ is nonempty. For a given $\mathbf{p}^+ \in P^+$, let $\mathbf{x}_i^+$ denotes a maximizer of the problem (2), $\mathbf{x}^+ = (\mathbf{x}_1^+,..,\mathbf{x}_n^+)$. Define $\pi_i^{*prod} = \pi_i^{prod}(\mathbf{p}^+, \mathbf{x}_i^*), \pi_i^{+prod} = \pi_i^{prod}(\mathbf{p}^+, \mathbf{x}_i^+)$, then the duality gap is given by

$$v - v^D = \sum_{i \in I} [\pi_i^{+prod} - \pi_i^{*prod}], \quad (4)$$

and, according to the approach developed in [20]-[22], represents the sum of generator lost profits corresponding to foregone opportunities to supply power in the amount corresponding to $\mathbf{x}^+$ at a price $\mathbf{p}^+ \in P^+$ by accepting dispatch $\mathbf{x}^*$. Since $\pi_i^{+prod} \geq \pi_i^{*prod}$, $\forall i \in I$, the duality gap is non-negative. If either of (1), (2), and (3) have multiple optimal points, the duality gap is independent of the choice of $\mathbf{x}^*$, $\mathbf{x}^+$, $\mathbf{p}^+ \in P^+$. The convex hull pricing method instructs to distribute the amount (4) to generators as uplift payments to ensure that no generator has an incentive to change its output given the market price $\mathbf{p}^+ \in P^+$ leaving aside issues related to market power. The following two interpretations are applicable to the set $P^+$. First, (4) entails that set of prices is chosen in a way to minimize the total uplift payment needed to support the centralized dispatch solution [20]. Second, the set $P^+$ can be viewed as the subdifferential of the convex function $v^D$, which is the convex hull of the total cost function $v$ viewed as a function of $\mathbf{d}$ [21], [22]. These insights justify the terms "minimum-uplift pricing" and "convex hull pricing" used to describe the method.

The case of price-sensitive demand bids is formulated analogously. For a consumer $j$, $j \in J$, $J = \{1,..,m\}$, $m < +\infty$, let $\mathbf{d}_j$, $\mathbf{d}_j \in R_{\geq 0}^T$, denote the vector of power demand. To consider the possibility of the demand-side non-convexities, let us introduce $J^{nc}$, $J^{nc} \subset J$, a subset of the consumers with non-convex demand and denote as $\mathbf{v}_j$, $j \in J^{nc}$, a vector of the corresponding discrete variables. Let a consumer $j$ bid the benefit function $B_j(\mathbf{y}_j)$ with $\mathbf{y}_j = (\mathbf{v}_j, \mathbf{d}_j)$ for $j \in J^{nc}$ and $\mathbf{y}_j = \mathbf{d}_j$, for $\forall j \in J \setminus J^{nc}$. For any fixed $\mathbf{v}_j$, $j \in J^{nc}$, $B_j(\mathbf{y}_j)$ viewed as a function of $\mathbf{d}_j$ is assumed to be the non-negative continuous concave function in its corresponding domain, which may depend on the value of $\mathbf{v}_j$. We also assume $B_j(\mathbf{y}_j)$ to be non-negative continuous concave function of $\mathbf{d}_j$, $\forall j \in J \setminus J^{nc}$. Let $Y_j$ denote the individual feasible set of a consumer $j$. Note that the set $Y_j$ may depend on the initial conditions referring to the past and/or the present market planning cycles. For simplicity, we do not explicitly indicate the initial conditions. The set $Y_j$ is assumed to be compact, $\forall j \in J$. Let $\mathbf{y} = (\mathbf{y}_1, ..., \mathbf{y}_m)$, the primal problem takes the form

$$U = \max_{\substack{\mathbf{x}, \mathbf{y} \\ \mathbf{x}_i \in X_i, \forall i \in I \\ \mathbf{y}_j \in Y_j, \forall j \in J \\ \sum_{i \in I} \mathbf{g}_i = \sum_{j \in J} \mathbf{d}_j}} \left[ \sum_{j \in J} B_j(\mathbf{y}_j) - \sum_{i \in I} C_i(\mathbf{x}_i) \right]. \quad (5)$$

The feasible set of (5) is denoted as $\Omega$ and is assumed to be nonempty and compact. Clearly, under the stated assumptions the maximum in the primal problem (5) exists. Let $\mathbf{x}^*, \mathbf{y}^*$ denote a maximizer of (5). The dual problem is formulated as

$$U^D = \min_{\mathbf{p} \in R^T} \max_{\substack{\mathbf{x}, \mathbf{y} \\ \mathbf{x}_i \in X_i, \forall i \in I \\ \mathbf{y}_j \in Y_j, \forall j \in J}} \left[ \mathbf{p}^T \left( \sum_{i \in I} \mathbf{g}_i - \sum_{j \in J} \mathbf{d}_j \right) + \sum_{j \in J} B_j(\mathbf{y}_j) - \sum_{i \in I} C_i(\mathbf{x}_i) \right] =$$

$$\min_{\mathbf{p} \in R^T} \left( \sum_{j \in J} \pi_j^{cons}(\mathbf{p}) + \sum_{i \in I} \pi_i^{prod}(\mathbf{p}) \right) \quad , (6)$$

with the consumer profit functions

$$\pi_j^{cons}(\mathbf{p}) = \max_{\mathbf{y}_j \in Y_j} \pi_j^{cons}(\mathbf{p}, \mathbf{y}_j), \quad \pi_j^{cons}(\mathbf{p}, \mathbf{y}_j) = B_j(\mathbf{y}_j) - \mathbf{p}^T \mathbf{d}_j. \quad (7)$$

Since (6) is unconstrained minimization problem of the convex function, its set of minimizers $P^+$ is given by the solutions to $\{0\} \in \sum_{j \in J} \partial \pi_j^{cons}(\mathbf{p}) + \sum_{i \in I} \partial \pi_i^{prod}(\mathbf{p})$. For a given $\mathbf{p}^+ \in P^+$, let us denote as $\mathbf{x}^+$ and $\mathbf{y}^+$ maximizers of (2) and (7) respectively and let $\pi_j^{*cons} = \pi_j^{cons}(\mathbf{p}^+, \mathbf{y}_j^*)$, $\pi_j^{+cons} = \pi_j^{cons}(\mathbf{p}^+, \mathbf{y}_j^+)$. The duality gap is given by

$$U^D - U = \sum_{j \in J} [\pi_j^{+cons} - \pi_j^{*cons}] + \sum_{i \in I} [\pi_i^{+prod} - \pi_i^{*prod}]. \quad (8)$$

The duality gap is non-negative and gives the total lost profit of consumers and producers due to the difference between the optimal values of output (consumption) volumes in the decentralized dispatch and these volumes set by the centralized dispatch for a given market price $\mathbf{p}^+ \in P^+$. In the convex hull pricing mechanism, the market participants are entitled to compensations of these lost profits to ensure equilibrium of the centralized dispatch solution. The problem of how to finance the calculated uplift payments is beyond the scope of the present paper, but we have the following comments on the issue. The primal problem solution may imply the cross-subsidization among the consumers: if it is optimal for the market to operate a unit, its output volume is set by the variable cost only without any reference to the start-up cost. If the unit's start-up cost is not fully compensated by the output sales at the market price and needs to be recovered through the uplift payment, in the general case the obligation to make this payment cannot be uniformly allocated among the consumers as no consumer should be charged the cost higher than indicated in its bid.

Thus, the allocation of the uplift charges among the market players is generally non-uniform. Also, if there is no fixed load segment of demand, then the combined profit of all the market participants equals $U$ - the value of the primal problem solution since the amount paid by consumers equals the amount received by producers. In this case, if the duality gap is present and no outside funds are supplied to the market, then it is impossible to compensate each market participant the corresponding lost profit obtained from (8). As a result, the calculated market uplift payments can be regarded as preliminary as these uplifts are to be further modified to ensure revenue adequacy. However, since the primal problem outcome ensures that the total consumer benefit is no lower than the total producer cost with some sufficiently low (high) cost (benefit) assigned to the output (consumption) volumes that are constrained by the initial conditions (as well as the fixed segment of demand), the lost profit (if any) associated with the opportunity to produce (consume) exactly these constrained volumes can be fully compensated by the market players without a need for the external funding. Thus, the market budget suffices to ensure the non-confiscatory pricing.

Clearly, if the primal and/or the dual problems have multiple solutions, then the duality gap is independent of the choice of both $\mathbf{x}^*$, $\mathbf{y}^*$ and $\mathbf{x}^+$, $\mathbf{y}^+$ as well as $\mathbf{p}^+ \in P^+$. To illustrate implications of the convex hull pricing we consider two examples. (For simplicity, throughout the present paper in case of uninode single-period market models we omit the vector notations and reference to the time period and assume that all the generating units are initially offline. Also, in this case the fixed cost corresponds to either start-up and/or no-load cost, for brevity we will refer to it as the start-up cost.)

*Example 1.*

Consider a uninode single-period power system with one consumer and one producer. The consumer benefit function is $B(d) = bd$, $0 \leq d \leq d^{\max}$, so that it submits one-step bid with a price $b$. The producer has cost function $C(x) = ag + wu$ with the start-up cost $w$, $w > 0$, and output $g$, $0 \leq g \leq g^{\max}$, so that its generating unit has zero minimum output limit and linear variable cost. The parameters are assumed to satisfy $a \geq 0$, $a + w/g^{\max} < b$, and $bd^{\max} < ad^{\max} + w$, which imply $d^{\max} < g^{\max}$. In this setting, the consumer is not able to compensate the total cost of power production for output in the range $0 \leq g \leq d^{\max}$, hence $u^* = g^* = d^* = 0$. The convex hull pricing results in the unique market price $p^+ = a + w/g^{\max}$ and zero uplift payment to the producer and the uplift of $(b - p^+)d^{\max}$ payable to the consumer (the source of funds for the payment is not discussed here). The convex hull pricing method states that the consumer has lost the opportunity to consume the power volume $d^{\max}$ by accepting the centralized dispatch and is entitled for the compensation. The mathematical reason for the nonzero consumer's uplift is the price dependence on $g^{\max}$, which is neither technologically feasible given the power balance constraint nor economically feasible for the stated benefit and cost functions. In addition, we observe that the consumer is not able to buy any power volume from the producer at all. Therefore, it had not deprived the consumer of any other market opportunities. As $g^{\max}$ increases, the market price tends to $a$ and uplift payment to the consumer grows, while no opportunities foregone by accepting the centralized dispatch outcome emerge for the market players.

*Example 2.*

Consider a uninode single-period power system with a single producer from the Example 1 and an (inverse) demand given by $p(d) = 2a(1 - d/d^{\max})$ with the consumption volume $d$, $0 \leq d \leq d^{\max}$. The latter corresponds to the consumer benefit function $B(d) = 2ad - ad^2/d^{\max}$. The parameters are assumed to satisfy $6w/a \leq d^{\max} < g^{\max}$, which ensures that the only solution of the primal problem is given by $u^* = 1$, $g^* = d^* = d^{\max}/2$ with $U = ad^{\max}/4 - w$. The marginal price equals $a$, which means that some segment of demand with prices equal or just above $a$ was cleared by the market. If $p(d)$ is an aggregate (inverse) demand of a pool of small consumers and no consumer is required to pay the cost above indicated in its bid, then some part of demand will not contribute to compensation of the generator's start-up cost and, therefore, the primal problem solution implies the cross-subsidization between the different segments of demand.

The solution for the dual problem yields the unique market price $p^+ = a + w/g^{\max}$, which gives $d^+ = [1 - w/(ag^{\max})]d^{\max}/2$ and two solutions for the generating unit status and output variables: the unit being online producing $g^+ = g^{\max}$ and the unit being offline with $g^+ = 0$. The uplift is given by

$$U^D - U = w\left(1 - \frac{d^{\max}}{2g^{\max}} + \frac{wd^{\max}}{a(2g^{\max})^2}\right), \quad (9)$$

which vanishes as $w \to 0$ since in the absence of the start-up cost the primal problem can be transformed into the convex optimization problem with no duality gap. We observe that both the price $p^+$ and the duality gap, which is equal to the total uplift (9), depend on $g^{\max}$. That seems counterintuitive since the technologically feasible output (consumption) volumes are given by the range $[0, d^{\max}]$ and the output volumes in the range $(d^{\max}, g^{\max}]$ exceed the maximum consumption volume $d^{\max}$ and thus are not technologically feasible. Also, the technologically feasible consumption volumes from the range $(d^*, d^{\max}]$ have prices below $a$ and, hence, are economically infeasible as this segment of demand is not able to cover even the marginal cost of output[2]. Hence, the technologically feasible output volumes in the range $(d^*, d^{\max}]$ are not economically feasible.

In these examples, the price $p^+$ dependence on $g^{\max}$, which exceeds technologically and/or economically feasible output volumes, raises a question if some opportunities, treated as foregone in the lost profit calculation (8), actually cannot be realized and, therefore, should not be taken into account in the lost profit calculation for a given value of the market price.

The examples above show that some opportunities that are considered as lost in the convex hull pricing mechanism are either technologically and/or economically infeasible, and hence the opportunities available to the market players in a centralized dispatch should be thoroughly examined. In the next section we identify the output (consumption) volumes that are economically and technologically feasible in the centralized dispatch and for each market player define the modified individual feasible set. These sets are further utilized in Section 4 in formulation of the proposed modified convex hull pricing method.

---

[2] Another economic constraint is $B(g) \geq ag + wu$, $0 \leq g \leq ud^{\max}$, $u \in \{0,1\}$, stating that the demand is able to pay for the generator total output cost. It is straightforward to check that this constraint excludes only some of the output volumes in the range $(d^*, d^{\max}]$, while the condition $p(g) \geq a$ discards all of them.

**3. The modified individual feasible sets of the market players**

To ensure the economic stability of the centralized dispatch outcome the lost profit compensations for the market players have to be introduced in the market. Since a market player not satisfied with the centralized dispatch outcome (including the corresponding pricing) may attempt to deviate from the centralized dispatch and consumption schedule in the real time, we need to specify pricing method for power sold/bought in the real-time market. We propose applying the following market principles:
1. A producer (consumer) is regarded as having an opportunity to supply (consume) some output (consumption) volume iff this volume could be scheduled in the centralized dispatch.
2. Any deviation from the generation or consumption schedule set by the centralized dispatch outside the certain error thresholds (except for that incurred due to participation in the power balancing) is economically discouraged.

The second principle implies that non-conformance of a market player with the centralized dispatch leads to deviation charges that exceed the expected profit from such actions. Thus, in order to have a desired output (consumption) schedule, the market player is encouraged to schedule it in the centralized dispatch. One way to achieve that is to sign the bilateral contracts with the other market players for this schedule.

To calculate the lost profit we need to identify a set of output (consumption) volumes that can be scheduled in the centralized dispatch and the corresponding profits at a given market price. We refer to the set of these volumes as technologically and economically feasible output (consumption) volumes for a producer (consumer).

In the examples 1 and 2 with one producer and one consumer it is natural to define the corresponding technologically and economically feasible output volumes for the producer as all output volumes that satisfy both the technological constraints (power balance constraints and individual market player constraints) and economic constraints (ability of demand to compensate the total production cost and ability of each segment of demand to pay the marginal cost of output), i.e. all output volumes satisfying $B(g) \geq C(u,g)$, $0 \leq g \leq ud^{\max}$, $u \in \{0,1\}$, and for nonzero output we also require $\partial_g B(g) \geq \partial_g C(u,g)$. The power balance constraints ensure that these are also the sets of technologically and economically feasible consumption volumes for the corresponding consumers. We observe that each output (consumption) volume belonging to the set for a given producer (consumer) can be obtained as a solution to the primal problem with some new values of $g^{\max}$ and $d^{\max}$ no higher than their original values.

To find technologically and economically feasible output (consumption) volumes in a more general setting let us consider a hypothetical set of bilateral contracts for power that a given market player may sign with the other market participant (both generators and consumers). We require the resulting contracted schedule to belong to the feasible set of the primal optimization problem, which implies that the nonzero output (consumption) volumes in the present market planning horizon constrained due to the market player individual operational constraints (such as fixed segment of demand, minimum output volumes of the units that were committed in the previous market planning cycle(s) and have to be online because of the minimum up time constraints) have to be fully contracted. In the considered bilateral contract framework, the market player pays the producers the full costs of contracted output according to their bids and collects payments from the consumers in the amounts equal to the values of their benefit functions at the respective contracted consumption volumes except for the abovementioned constrained output (consumption) volumes with the sunk cost (benefit), which can be formally assigned any arbitrary low (high) value.

In this hypothetical setting, for every time period of the given market planning cycle only the market participants with contracts are allowed to participate in the centralized dispatch optimization problem with their original bids and volumes not exceeding the contracted output (consumption) volumes, while the market player in question submits its original bid with volume not exceeding the netted contracted volume for each time period. Thus, these contracted volumes are not mandatory for inclusion into the dispatch (they can be viewed as the financial contracts accounted for at the financial settlement stage only), while uncontracted volumes do not participate in the dispatch. We will use this theoretical bilateral contract construction to identify the output (consumption) volumes that can be supplied (consumed) in the centralized dispatch, i.e. technologically and economically feasible output (consumption) volumes.

We emphasize that in the general case this approach is more restrictive than merely adding the requirement that the consumers are able to pay the producers the total output cost with the proper treatment of the sunk cost (benefit) incurring in the present market planning horizon due to the initial condition for the market players.

Consider the case when demand is fixed and given by the vector $\mathbf{d}$. We say that a generator $i_0$ has an opportunity to supply an output volume $\overline{\mathbf{g}}_{i_0}$ if the generator can sign the bilateral contracts with the other market players (both producers and consumers) with $i$-th generator having an output volume $\overline{\mathbf{g}}_i$ and the corresponding status $\overline{\mathbf{u}}_i$, $\overline{\mathbf{x}}_i \in X_i$, so that the load is fully contracted and the resulting regime can be realized as the centralized dispatch. This is formalized as the following condition for $\overline{\mathbf{x}} = (\overline{\mathbf{x}}_1,...,\overline{\mathbf{x}}_n)$:

$$\overline{\mathbf{x}} \in \arg \min_{\substack{\mathbf{x}, \\ \mathbf{x}_i \in X_i, \forall i \in I \\ \mathbf{x} \leq \overline{\mathbf{x}} \\ \sum_{i \in I} \mathbf{g}_i = \mathbf{d}}} \sum_{i \in I} C_i(\mathbf{x}_i), \quad (10)$$

where $\mathbf{x} \leq \overline{\mathbf{x}}$ denotes a set of the inequality constraints $u_i^t \leq \overline{u}_i^t$, $g_i^t \leq \overline{g}_i^t$, $\forall i \in I$, $\forall t \in T$. From (10) it follows that each $\overline{\mathbf{x}}$ belongs to the feasible set of the primal problem: $\overline{\mathbf{x}} \in \Omega$. Let $\overline{\Omega}$ be the set of all $\overline{\mathbf{x}}$ that satisfy (10). Since, in addition to supplying output $\overline{\mathbf{g}}_{i_0}$, the producer $i_0$ merely acts as an intermediary (retaining the market surplus with the proper treatment of the sunk cost (benefit)), the set $\overline{\Omega}$ is independent of the choice of $i_0$. For a producer $i'$, we denote as $\overline{\Omega}_{i'}^{prod}$ a set of all $\overline{\mathbf{x}}_{i'} = (\overline{\mathbf{u}}_{i'}, \overline{\mathbf{g}}_{i'})$ such that $\mathbf{x} \in \overline{\Omega}$ for some $\mathbf{x}_i$, $\forall i \in I$, $i \neq i'$. Equivalently, the set $\overline{\Omega}_{i'}^{prod}$ can be defined as a projection of $\overline{\Omega}$ on $X_{i'}$. The set of these $\overline{\mathbf{g}}_{i'}$ is referred to as the set of the technologically and economically feasible output volumes for a producer $i'$. Generally, $\overline{\Omega}_i^{prod} \subset X_i$, $\overline{\Omega} \subset \times_{i \in I} \overline{\Omega}_i^{prod}$, but $\overline{\Omega} \neq \times_{i \in I} \overline{\Omega}_i^{prod}$. Clearly, the original centralized dispatch outcome $\mathbf{x}^*$ can be realized through a set of bilateral contracts, since it satisfies (10). Thus, we have $\mathbf{x}^* \in \overline{\Omega} \subset \Omega$ and $\mathbf{x}_i^* \in \overline{\Omega}_i^{prod}$, $\forall i \in I$. We also note that the resulting sets $\overline{\Omega}$ and $\overline{\Omega}_i^{prod}$, $\forall i \in I$, generally depend on the initial conditions for a given market planning horizon.

We note that (10) can be reformulated to exclude $\overline{\mathbf{x}}$ from the constraint set of the optimization problem. Clearly, $\overline{\mathbf{x}}$ satisfies (10) iff

$$\overline{\mathbf{x}} \in \arg \min_{\substack{\mathbf{x}, \\ \mathbf{x}_i \in X_i, \forall i \in I \\ \mathbf{x} \leq \mathbf{x}^\circ \\ \sum_{i \in I} \mathbf{g}_i = \mathbf{d}}} \sum_{i \in I} C_i(\mathbf{x}_i)$$

for some $\mathbf{x}^\circ$. Thus, the condition (10) can be interpreted in the following intuitive way. For definiteness, consider the case when each time period of a given market planning cycle is given by one hour. For each generator in the power system, let us take any set of hourly output volumes allowed by the generator individual operational constraints as new hourly maximum output limits of the generator (without changing any other technical parameters of the generating units) and consider the new primal problem. If the resulting centralized dispatch problem is feasible, it will produce some new hourly output volumes for each generator. Clearly, these output volumes with the appropriate values of the binary status variables satisfy (10). As we simultaneously vary the set of new hourly maximum output limits of all the generators in the power system and find the new centralized dispatch outputs, for each generator we obtain a set of technologically and economically feasible output volumes (with $\overline{\Omega}_i^{prod}$ being the corresponding status-output set). In this setting, the hypothetical bilateral contracts are utilized to provide an economic framework in which the market players participate in the centralized dispatch market not only with their actual values of the hourly maximum output limits but also with the reduced ones.

Another interpretation of (10) is that starting with value of the objective function $\sum_{i \in I} C_i(\mathbf{x}_i)$ at a point $\overline{\mathbf{x}}$, this value cannot be increased by reducing the output volumes of the generators (without increasing output of some generator) staying within a feasible set of (1).

Let $C_i^{sunk}$ denote the sunk cost of a generator $i$ in the present market planning cycle[3]. For each unit $i$, introduce a set $\Psi_i^{prod}$:

$$\Psi_i^{prod} = \{\mathbf{x}_i \mid \mathbf{x}_i \in X_i, C_i(\mathbf{x}_i) = C_i^{sunk}\}. \quad (11)$$

Clearly, the set $\Psi_i^{prod}$ satisfies $\Psi_i^{prod} \subset X_i$ for the given initial condition for a unit $i$ and reflects the generator's option not to participate in the market[4] except for the cases of units being online due to previous commitment that cannot be altered, must-

---

[3] For the case of non-negative pricing, it is natural to define the sunk cost as $C_i^{sunk} = \min_{\mathbf{x}_i \in X_i} C_i(\mathbf{x}_i)$.

[4] Alternative (non-equivalent) definition of $\Psi_i^{prod}$, which also ensures non-confiscatory pricing, is a subset of (11) defined as $\Psi_i^{prod} = \bigcup_{\substack{\mathbf{p}, \\ \mathbf{p} \in R_{\geq 0}^T, \\ \|\mathbf{p}\|=1}} \arg \min_{\substack{\mathbf{x}_i, \\ \mathbf{x}_i \in X_i \\ C_i(\mathbf{x}_i) = C_i^{sunk}}} \mathbf{p}^T \mathbf{g}_i$.

run statuses, etc. We propose defining the modified individual feasible set for a generator $i$ as follows:

$$\widetilde{X}_i(\boldsymbol{\varepsilon}_i) = \Psi_i^{prod} \cup \bigcup_{\overline{\mathbf{x}}_i \in \overline{\Omega}_i^{prod}} \Delta_i^{prod}(\overline{\mathbf{x}}_i, \boldsymbol{\varepsilon}_i), \forall i \in I, \quad (12)$$

where $\boldsymbol{\varepsilon}_i \in R_{>0}^T$ (i.e., $\boldsymbol{\varepsilon}_i = (\varepsilon_i^1, ..., \varepsilon_i^T)$, $\varepsilon_i^t > 0$, $\forall t \in \{1, .., T\}$), and

$$\Delta_i^{prod}(\overline{\mathbf{x}}_i, \boldsymbol{\varepsilon}_i) = \{\mathbf{x}_i \mid \mathbf{x}_i \in X_i, \left| g_i^t - \overline{g}_i^t \right| \leq \varepsilon_i^t, \forall t \in \{1, .., T\}\}, \forall i \in I. \quad (13)$$

Utilization of the set $\Psi_i^{prod}$ in (12) ensures the non-confiscatory pricing for the unit $i$ as its profit obtained from the dual problem solution for a market price $\mathbf{p}$ is no lower than $\pi_i^{prod}(\mathbf{p}, \mathbf{x}_i)$, $\forall \mathbf{x}_i \in \Psi_i^{prod}$. If for each time interval of the given market planning period the market price in non-negative, then the loss (if any) of the generator $i$ is limited by the sunk cost $C_i^{sunk}$. If minimum of a generator $i$ cost is attained only when it is offline in all the time periods of a given market planning cycle and this is allowed by the unit's operational constraints for the given initial conditions, then $C_i^{sunk} = 0$ and (11) implies that $\Psi_i^{prod}$ contains only the element $\mathbf{u}_i = \mathbf{g}_i = 0$. In this case, definition (12) of $\widetilde{X}_i(\boldsymbol{\varepsilon}_i)$ guarantees non-negativity of the unit $i$ profit.

The inclusion of the set $\Delta_i^{prod}(\overline{\mathbf{x}}_i, \boldsymbol{\varepsilon}_i)$, $\boldsymbol{\varepsilon}_i \in R_{>0}^T$, for each element $\overline{\mathbf{x}}_i \in \overline{\Omega}_i^{prod}$, is needed to indicate in the dual problem whether at a given price the generator is willing to supply some more/less power than $\overline{g}_i^t$ belonging to its individual feasible set $X_i$. Thus, for each $\overline{\mathbf{g}}_i$ with $\overline{\mathbf{x}}_i \in \overline{\Omega}_i^{prod}$, the sets $\widetilde{X}_i(\boldsymbol{\varepsilon}_i)$ and $X_i$ have identical output volumes in some closed $\boldsymbol{\varepsilon}_i$-neighborhood of $\overline{\mathbf{g}}_i$. By construction, we have $\widetilde{X}_i(\boldsymbol{\varepsilon}_i) \subset X_i$.

Now we turn to the case when the price-sensitive consumer bids are present in the power market. We say that a generator $i_0$ has an opportunity to supply an output volume $\overline{\mathbf{g}}_{i_0}$ if the generator can sign the bilateral contracts with the other market players so that the generators have status-output values $\overline{\mathbf{x}}$ with $\overline{\mathbf{x}}_i \in X_i$, $\forall i \in I$, and the consumers have consumption specified by $\overline{\mathbf{y}} = (\overline{\mathbf{y}}_1, ..., \overline{\mathbf{y}}_m)$, $\overline{\mathbf{y}}_j \in Y_j$, $\forall j \in J$, satisfying the following condition:

$$\{\overline{\mathbf{x}}, \overline{\mathbf{y}}\} \in \arg \max_{\substack{\mathbf{x}, \mathbf{y} \\ \mathbf{x}_i \in X_i, \forall i \in I \\ \mathbf{y}_j \in Y_j, \forall j \in J \\ \mathbf{x} \leq \overline{\mathbf{x}}, \mathbf{y} \leq \overline{\mathbf{y}} \\ \sum_{i \in I} \mathbf{g}_i = \sum_{j \in J} \mathbf{d}_j}} \left[ \sum_{j \in J} B_j(\mathbf{y}_j) - \sum_{i \in I} C_i(\mathbf{x}_i) \right]. \quad (14)$$

The set of technologically and economically feasible output (consumption) volumes for the generator $i$ (consumer $j$) is given by all $\overline{\mathbf{g}}_i$ ($\overline{\mathbf{d}}_j$) that satisfy (14). Clearly, we have $\{\overline{\mathbf{x}}, \overline{\mathbf{y}}\} \in \Omega$. Likewise, let $\overline{\Omega}$ denote the set of such $\{\overline{\mathbf{x}}, \overline{\mathbf{y}}\}$. Instead of using generator $i_0$ to assemble the possible centralized dispatch outcome, we could use any other producer (consumer) for that purpose with the resulting set $\overline{\Omega}$ being independent of this choice. Analogously, for each producer $i$ and consumer $j$ we define the sets $\overline{\Omega}_i^{prod}$ and $\overline{\Omega}_j^{cons}$ that are projections of $\overline{\Omega}$ on the set $X_i$ and $Y_j$, respectively. For any $i$, $j$ we have $\mathbf{x}_i^* \in \overline{\Omega}_i^{prod} \subset X_i$ and $\mathbf{y}_j^* \in \overline{\Omega}_j^{cons} \subset Y_j$. Also, in the general case $\overline{\Omega} \subset \left(\underset{i \in I}{\times} \overline{\Omega}_i^{prod}\right) \times \left(\underset{j \in J}{\times} \overline{\Omega}_j^{cons}\right)$, but $\overline{\Omega} \neq \left(\underset{i \in I}{\times} \overline{\Omega}_i^{prod}\right) \times \left(\underset{j \in J}{\times} \overline{\Omega}_j^{cons}\right)$. We also have $\{\mathbf{x}^*, \mathbf{y}^*\} \in \overline{\Omega} \subset \Omega$. As in the fixed demand case, the resulting sets $\overline{\Omega}$, $\overline{\Omega}_i^{prod}$, and $\overline{\Omega}_j^{cons}$ generally depend on the initial conditions for a given market planning cycle.

As in the fixed load case, (14) can be equivalently formulated without use of $\overline{\mathbf{x}}$, $\overline{\mathbf{y}}$ in the feasible set of (14). Likewise, (14) implies that a set of technologically and economically feasible output (consumption) volumes for a generator $i$ (consumer $j$) is given by the solutions to the centralized dispatch problem with new hourly maximum output and consumption limits of all the market players not exceeding the corresponding limits used in the primal problem. Also, (14) implies that starting with value of the total market utility function at a point $\{\overline{\mathbf{x}}, \overline{\mathbf{y}}\}$, this value cannot be increased by reducing the hourly output (consumption) volumes of the generators (consumer) without increasing the hourly output (consumption) of some generator (consumer) provided that we stay within a feasible set of (5).

Let $B_j^{sunk}$ denote the sunk benefit of a consumer $j$ in the present market planning cycle. For a consumer $j$, the set $\Psi_j^{cons}$ is introduced to ensure the non-confiscatory pricing for power except for the sunk benefit $B_j^{sunk}$ and is given by

$$\Psi_j^{cons} = \{\mathbf{y}_j \mid \mathbf{y}_j \in Y_j, B_j(\mathbf{y}_j) = B_j^{sunk}\}. \quad (15)$$

Similarly, we define the modified individual feasible set for a market player in case of price-sensitive consumer bids. The set for a generator $i$ is defined by (12) with some $\boldsymbol{\varepsilon}_i \in R_{>0}^T$, while $\widetilde{Y}_j(\boldsymbol{\rho}_j)$ - the modified individual feasible set for a consumer $j$ is given by[5]

$$\widetilde{Y}_j(\boldsymbol{\rho}_j) = \Psi_j^{cons} \cup \bigcup_{\overline{\mathbf{y}}_j \in \overline{\Omega}_j^{cons}} \Delta_j^{cons}(\overline{\mathbf{y}}_j, \boldsymbol{\rho}_j),$$

where $\Delta_j^{cons}(\overline{\mathbf{y}}_j, \boldsymbol{\rho}_j) = \{\mathbf{y}_j \mid \mathbf{y}_j \in Y_j, |d_j^t - \overline{d}_j^t| \leq \rho_j^t, \forall t \in \{1,..,T\}\}$ with some $\boldsymbol{\rho}_j \in R_{>0}^T$, $\forall j \in J$. The expression (15) for $\Psi_j^{cons}$ ensures the non-confiscatory pricing for a consumer $j$ because its profit obtained from the dual problem solution for a market price $\mathbf{p}$ is no lower than $\pi_j^{cons}(\mathbf{p}, \mathbf{y}_j)$, $\forall \mathbf{y}_j \in \Psi_j^{cons}$. In the simple case, when a consumer $j$ with no demand-side non-convexity has price-insensitive consumption volume $\mathbf{d}_j^{min}$ and consumes some additional power only if the prices are sufficiently low, the set $\Psi_j^{cons}$ contains only the element $\mathbf{d}_j = \mathbf{d}_j^{min}$ and $B_j^{sunk} = B_j(\mathbf{d}_j^{min})$. We note that $\widetilde{Y}_j(\boldsymbol{\rho}_j)$ may not be compact.

We also note that since all $\overline{\mathbf{x}}$ in the fixed load case and all $\{\overline{\mathbf{x}}, \overline{\mathbf{y}}\}$ in the case with price-sensitive consumer bids are primal feasible, the right-hand sides of (10) and (14) define the respective self-correspondence $\Omega \to\to \Omega$ with a range $\overline{\Omega}$. Thus, (10) and (14) imply that $\overline{\Omega}$ is a set of fixed points of that correspondence.

---

[5] Likewise, to ensure the non-confiscatory pricing the set $\Psi_j^{cons}$ can be alternatively defined as $\Psi_j^{cons} = \bigcup_{\substack{\mathbf{p}, \\ \mathbf{p} \in R_{\geq 0}^T, \\ \|\mathbf{p}\|=1}} \underset{\substack{\mathbf{y}_j, \\ \mathbf{y}_j \in Y_j, \\ B_j(\mathbf{y}_j) = B_j^{sunk}}}{\arg\min} \mathbf{p}^T \mathbf{d}_j$, which is a subset of (15).

## 4. The modified convex hull pricing for systems with the fixed load

Let us define the modified primal problem
$$\tilde{v}(\boldsymbol{\varepsilon}) = \min_{\substack{\mathbf{x}, \\ \mathbf{x}_i \in \tilde{X}_i(\boldsymbol{\varepsilon}_i), \forall i \in I \\ \sum_{i \in I} \mathbf{g}_i = \mathbf{d}}} \sum_{i \in I} C_i(\mathbf{x}_i) \quad (16)$$

with some $\boldsymbol{\varepsilon} = (\boldsymbol{\varepsilon}_1,..,\boldsymbol{\varepsilon}_n)$, $\boldsymbol{\varepsilon}_i \in R^T_{>0}$, $\boldsymbol{\varepsilon} \in R^{nT}_{>0}$. Let $\tilde{\Omega}(\boldsymbol{\varepsilon})$ denote the feasible set of the modified primal problem (16). Since $\mathbf{x}^* \in \tilde{\Omega}(\boldsymbol{\varepsilon})$ for any $\mathbf{x}^*$ - an optimal point for the primal problem (1) and $\tilde{\Omega}(\boldsymbol{\varepsilon}) \subset \Omega$, we conclude that $v = \tilde{v}(\boldsymbol{\varepsilon})$. Consequently, both the primal and the modified primal problems have the identical sets of minimizers. The dual of the modified primal problem is
$$\tilde{v}^D(\boldsymbol{\varepsilon}) = \max_{\mathbf{p} \in R^T} \left( \mathbf{p}^T \mathbf{d} - \sum_{i \in I} \tilde{\pi}_i^{prod}(\mathbf{p}, \boldsymbol{\varepsilon}_i) \right), \quad (17)$$

with
$$\tilde{\pi}_i^{prod}(\mathbf{p}, \boldsymbol{\varepsilon}_i) = \sup_{\mathbf{x}_i \in \tilde{X}_i(\boldsymbol{\varepsilon}_i)} \pi_i^{prod}(\mathbf{p}, \mathbf{x}_i). \quad (18)$$

We note that $\tilde{X}_i(\boldsymbol{\varepsilon}_i)$ is bounded and $\pi_i^{prod}(\mathbf{p}, \mathbf{x}_i)$ is the continuous function of $\mathbf{g}_i$ on $cl[\tilde{X}_i(\boldsymbol{\varepsilon}_i)]$ for each fixed $\mathbf{u}_i$, where $cl$ denotes closure of a set. (Note that since $X_i$ is closed and $\tilde{X}_i(\boldsymbol{\varepsilon}_i) \subset X_i$, we have $cl[\tilde{X}_i(\boldsymbol{\varepsilon}_i)] \subset X_i$.) Therefore:
$$\tilde{\pi}_i^{prod}(\mathbf{p}, \boldsymbol{\varepsilon}_i) = \max_{\mathbf{x}_i \in cl[\tilde{X}_i(\boldsymbol{\varepsilon}_i)]} \pi_i^{prod}(\mathbf{p}, \mathbf{x}_i), \quad \forall i \in I. \quad (19)$$

Since each $\tilde{\pi}_i^{prod}(\mathbf{p}, \boldsymbol{\varepsilon}_i)$ is a point-wise maximum of the function linear in $\mathbf{p}$, it is convex in $\mathbf{p}$ with well-defined subdifferential with respect to $\mathbf{p}$ denoted as $\partial \tilde{\pi}_i^{prod}(\mathbf{p}, \boldsymbol{\varepsilon}_i)$. It is straightforward to verify that a set of maximizers for (17), which we denote as $\tilde{P}^+(\boldsymbol{\varepsilon})$, is nonempty. Clearly, $\tilde{P}^+(\boldsymbol{\varepsilon})$ is a set of prices that solve $\mathbf{d} \in \sum_{i \in I} \partial \tilde{\pi}_i^{prod}(\mathbf{p}, \boldsymbol{\varepsilon}_i)$. For $\tilde{\mathbf{p}}^+ \in \tilde{P}^+(\boldsymbol{\varepsilon})$, let $\tilde{\mathbf{x}}_i^+$ denote a maximizer for the problem (19), and let us define $\tilde{\pi}_i^{*prod}(\boldsymbol{\varepsilon}) = \pi_i^{prod}(\tilde{\mathbf{p}}^+, \mathbf{x}_i^*)$, $\tilde{\pi}_i^{+prod}(\boldsymbol{\varepsilon}) = \pi_i^{prod}(\tilde{\mathbf{p}}^+, \tilde{\mathbf{x}}_i^+)$. The duality gap is given by
$$\tilde{v}(\boldsymbol{\varepsilon}) - \tilde{v}^D(\boldsymbol{\varepsilon}) = \sum_{i \in I} [\tilde{\pi}_i^{+prod}(\boldsymbol{\varepsilon}) - \tilde{\pi}_i^{*prod}(\boldsymbol{\varepsilon})].$$

Relation $\tilde{X}_i(\boldsymbol{\varepsilon}_i) \subset X_i$, $\forall i \in I$, implies $v^D \leq \tilde{v}^D(\boldsymbol{\varepsilon})$. Hence, $v^D \leq \tilde{v}^D(\boldsymbol{\varepsilon}) \leq \tilde{v}(\boldsymbol{\varepsilon}) = v$, which entails the following relation between the duality gaps of the original and the modified optimization problems:
$$0 \leq \tilde{v}(\boldsymbol{\varepsilon}) - \tilde{v}^D(\boldsymbol{\varepsilon}) \leq v - v^D, \quad \forall \boldsymbol{\varepsilon} \in R^{nT}_{>0}.$$

Therefore, the total uplift needed to support the centralized dispatch solution $\forall \tilde{\mathbf{p}}^+ \in \tilde{P}^+(\boldsymbol{\varepsilon})$ in the modified optimization problem (16) is no higher than that for the original problem (1), $\forall \mathbf{p}^+ \in P^+$. We propose computing the set of market prices and the individual uplifts $[\tilde{\pi}_i^{+prod}(\boldsymbol{\varepsilon}) - \tilde{\pi}_i^{*prod}(\boldsymbol{\varepsilon})]$ in the limit as $\boldsymbol{\varepsilon} \to +0$. We note that if there is an obligation to compensate the sunk cost (for example, if this cost originates from the commitment in the previous market planning period), then such an obligation is not a part of the considered total uplift payment [20].

A price $\mathbf{p}$ is said to support the primal problem solution $\mathbf{x}^*$ in the decentralized dispatch problem if

$$\mathbf{x}_i^* \in \arg\max_{\mathbf{x}_i \in X_i}\left[\mathbf{p}^T\mathbf{g}_i - C_i(\mathbf{x}_i)\right], \forall i \in I. \quad (20)$$

It is straightforward to verify that a price that supports the primal problem solution in the decentralized dispatch problem exists iff the duality gap is zero. In this case, the set of prices that support solution $\mathbf{x}^*$ is the same for all $\mathbf{x}^*$ (if the primal problem has multiple solutions) and is identical to a set of maximizers for the dual optimization problem. The substitution of $X_i$ by $cl[\widetilde{X}_i(\boldsymbol{\varepsilon}_i)]$ in (20) yields the definition of a price that supports the solution $\mathbf{x}^*$ in the modified decentralized dispatch problem.

Let us consider the case of uninode multi-period power system without non-convexities: no fixed costs, zero minimum output limits of all the units, no minimum up/down time, etc. Since in this case the primal problem can be transformed into the convex optimization problem with zero duality gap, the total uplift payment equals zero in both methods. The sets of prices that support solution $\mathbf{x}_i^*$ in the decentralized dispatch problem and the modified decentralized dispatch problem are identical as implied by the following proposition.

*Proposition 1:* If for some $i$ after optimization over the binary variables $\mathbf{u}_i$ the function $C_i(\mathbf{x}_i)$ becomes convex function of $\mathbf{g}_i$, then a price $\mathbf{p}$ supports solution $\mathbf{x}_i^*$ in the decentralized dispatch problem iff it supports solution $\mathbf{x}_i^*$ in the modified decentralized dispatch problem, $\forall \boldsymbol{\varepsilon}_i \in R_{\geq 0}^T$:

$$\mathbf{x}_i^* \in \arg\max_{\mathbf{x}_i \in X_i}\left[\mathbf{p}^T\mathbf{g}_i - C_i(\mathbf{x}_i)\right] \Leftrightarrow \mathbf{x}_i^* \in \arg\max_{\mathbf{x}_i \in cl[\widetilde{X}_i(\boldsymbol{\varepsilon}_i)]}\left[\mathbf{p}^T\mathbf{g}_i - C_i(\mathbf{x}_i)\right]. \quad (21)$$

*Proof.* Clearly, since $\mathbf{x}_i^* \in \widetilde{X}_i(\boldsymbol{\varepsilon}_i) \subset cl[\widetilde{X}_i(\boldsymbol{\varepsilon}_i)] \subset X_i$, if $\mathbf{p}$ supports $\mathbf{x}_i^*$ in the decentralized dispatch problem, then it supports $\mathbf{x}_i^*$ in the modified decentralized dispatch problem. To prove the converse, note that since $\mathbf{x}_i^* \in \widetilde{X}_i(\boldsymbol{\varepsilon}_i)$, we have

$$\mathbf{x}_i^* \in \arg\max_{\mathbf{x}_i \in cl[\widetilde{X}_i(\boldsymbol{\varepsilon}_i)]}\left[\mathbf{p}^T\mathbf{g}_i - C_i(\mathbf{x}_i)\right] \Rightarrow \mathbf{x}_i^* \in \arg\max_{\mathbf{x}_i \in \widetilde{X}_i(\boldsymbol{\varepsilon}_i)}\left[\mathbf{p}^T\mathbf{g}_i - C_i(\mathbf{x}_i)\right]. \quad (22)$$

For each fixed $\mathbf{g}_i \in R_{\geq 0}^T$, let us define

$$c_i(\mathbf{g}_i) = \min_{\substack{\mathbf{u}_i, \\ \mathbf{x}_i \in X_i}} C_i(\mathbf{x}_i). \quad (23)$$

(We adopt the convention that the value function equals $+\infty$ if the minimization problem is infeasible.) By assumptions, $c_i(\mathbf{g}_i)$ and $dom[c_i]$ are the convex function and the convex set, respectively. Likewise, for each fixed $\mathbf{g}_i \in R_{\geq 0}^T$, define

$$\widetilde{c}_i(\mathbf{g}_i) = \min_{\substack{\mathbf{u}_i, \\ \mathbf{x}_i \in \widetilde{X}_i(\boldsymbol{\varepsilon}_i)}} C_i(\mathbf{x}_i).$$

We note that $dom[\widetilde{c}_i](\boldsymbol{\varepsilon}_i)$ generally depends on $\boldsymbol{\varepsilon}_i$. Also, $\widetilde{X}_i(\boldsymbol{\varepsilon}_i) \subset X_i$ implies that if for a given $\mathbf{g}_i$ we have $\mathbf{x}_i \in \widetilde{X}_i(\boldsymbol{\varepsilon}_i)$ for some $\mathbf{u}_i$, then $\mathbf{x}_i \in X_i$. Therefore, $dom[\widetilde{c}_i](\boldsymbol{\varepsilon}_i) \subset dom[c_i]$.

We also have $c_i(\mathbf{g}_i) = \tilde{c}_i(\mathbf{g}_i)$, $\forall \mathbf{g}_i \in dom[\tilde{c}_i](\boldsymbol{\varepsilon}_i)$. Indeed, from (13) it follows that for a given $\mathbf{g}_i$ we have

$$\{\mathbf{u}_i \mid (\mathbf{u}_i, \mathbf{g}_i) \in \bigcup_{\overline{\mathbf{x}}_i \in \overline{\Omega}_i^{prod}} \Delta_i^{prod}(\overline{\mathbf{x}}_i, \boldsymbol{\varepsilon}_i)\} = \{\mathbf{u}_i \mid (\mathbf{u}_i, \mathbf{g}_i) \in X_i\}.$$

Consequently, if $(\mathbf{u}_i, \mathbf{g}_i) \in \bigcup_{\overline{\mathbf{x}}_i \in \overline{\Omega}_i^{prod}} \Delta_i^{prod}(\overline{\mathbf{x}}_i, \boldsymbol{\varepsilon}_i)$ for some $\mathbf{u}_i$, then $c_i(\mathbf{g}_i) = \tilde{c}_i(\mathbf{g}_i)$. Also, if $(\mathbf{u}_i, \mathbf{g}_i) \in \Psi_i^{prod}$, then $c_i(\mathbf{g}_i) = \tilde{c}_i(\mathbf{g}_i) = C_i^{sunk}$.

Moreover, (13) entails that $\forall \overline{\mathbf{g}}_i$ we have $\overline{\mathbf{g}}_i \in dom[\tilde{c}_i](\boldsymbol{\varepsilon}_i)$ and $dom[\tilde{c}_i](\boldsymbol{\varepsilon}_i) = dom[c_i]$ in the closed $\boldsymbol{\varepsilon}_i$-neighborhood of $\overline{\mathbf{g}}_i$ (note that this neighborhood may not be a subset of $dom[c_i]$). Since $\mathbf{x}_i^* \in \overline{\Omega}_i^{prod}$, this implies that the sets $dom[\tilde{c}_i](\boldsymbol{\varepsilon}_i)$ and $dom[c_i]$ are identical in the closed $\boldsymbol{\varepsilon}_i$-neighborhood of $\mathbf{g}_i^*$. Let $\mathbf{p}$ support $\mathbf{x}_i^*$ in the modified decentralized dispatch problem (19) for generator $i$. Let us optimize over $\mathbf{u}_i$ in (22), then the concave function $[\mathbf{p}^T \mathbf{g}_i - c_i(\mathbf{g}_i)]$ has a local maximum at $\mathbf{g}_i^*$ on a convex set $dom[c_i]$. Therefore, it has a global maximum at $\mathbf{g}_i^*$ on $dom[c_i]$, which together with (23) entails that $\mathbf{p}$ supports $\mathbf{x}_i^*$ in the decentralized dispatch problem (2) and, hence, (21) holds. The proposition is proved.

We conclude that for a convex power system both the convex hull pricing and the modified convex hull pricing methods result in identical sets of prices, which are the marginal prices. Also, Proposition 1 holds if $\mathbf{x}_i^*$ is replaced by $\forall \overline{\mathbf{x}}_i \in \overline{\Omega}_i^{prod}$.

## 5. The modified convex hull pricing for power systems with price-sensitive load

Now, let us consider the case when price-sensitive consumer bids are present as well. The modified primal problem is defined as:

$$\tilde{U}(\boldsymbol{\varepsilon}, \boldsymbol{\rho}) = \max_{\substack{\mathbf{x}, \mathbf{y} \\ \mathbf{x}_i \in \tilde{X}_i(\boldsymbol{\varepsilon}_i), \forall i \in I \\ \mathbf{y}_j \in \tilde{Y}_j(\boldsymbol{\rho}_j), \forall j \in J \\ \sum_{i \in I} \mathbf{g}_i = \sum_{j \in J} \mathbf{d}_j}} \left[ \sum_{j \in J} B_j(\mathbf{y}_j) - \sum_{i \in I} C_i(\mathbf{x}_i) \right] \quad (24)$$

with some $\boldsymbol{\varepsilon} \in R_{>0}^{nT}$, and $\boldsymbol{\rho} = (\boldsymbol{\rho}_1, ..., \boldsymbol{\rho}_m)$, $\boldsymbol{\rho}_j \in R_{>0}^T$, $\boldsymbol{\rho} \in R_{>0}^{mT}$. Let $\tilde{\Omega}(\boldsymbol{\varepsilon}, \boldsymbol{\rho})$ denote the feasible set of the modified primal problem (24). As $\tilde{\Omega}(\boldsymbol{\varepsilon}, \boldsymbol{\rho}) \subset \Omega$ and $(\mathbf{x}^*, \mathbf{y}^*) \in \tilde{\Omega}(\boldsymbol{\varepsilon}, \boldsymbol{\rho})$ for any $(\mathbf{x}^*, \mathbf{y}^*)$, which is maximizer of the primal problem (5), we deduce that the sets of maximizers for both the primal and the modified primal problems are identical and, therefore, $U = \tilde{U}(\boldsymbol{\varepsilon}, \boldsymbol{\rho})$. The dual of the modified primal problem is formulated as follows:

$$\tilde{U}^D(\boldsymbol{\varepsilon}, \boldsymbol{\rho}) = \min_{\mathbf{p} \in R^T} \left( \sum_{j \in J} \tilde{\pi}_j^{cons}(\mathbf{p}, \boldsymbol{\rho}_j) + \sum_{i \in I} \tilde{\pi}_i^{prod}(\mathbf{p}, \boldsymbol{\varepsilon}_i) \right), \quad (25)$$

with

$$\tilde{\pi}_j^{cons}(\mathbf{p}, \boldsymbol{\rho}_j) = \sup_{\mathbf{y}_j \in \tilde{Y}_j(\boldsymbol{\rho}_j)} \pi_j^{cons}(\mathbf{p}, \mathbf{y}_j)$$

and $\tilde{\pi}_i^{prod}(\mathbf{p}, \boldsymbol{\varepsilon}_i)$ given by (18). The functions $\pi_j^{cons}(\mathbf{p}, \mathbf{y}_j)$ and $\pi_i^{prod}(\mathbf{p}, \mathbf{x}_i)$ are continuous in $\mathbf{d}_j$ and $\mathbf{g}_i$ on $cl[\tilde{Y}_j(\boldsymbol{\rho}_j)]$ and $cl[\tilde{X}_i(\boldsymbol{\varepsilon}_i)]$ for each fixed $\mathbf{v}_j$ (if any) and $\mathbf{u}_i$, respectively. The sets $\tilde{Y}_j(\boldsymbol{\rho}_j)$, $\forall j \in J$, $\tilde{X}_i(\boldsymbol{\varepsilon}_i)$, $\forall i \in I$, are bounded, and we have

$$\tilde{\pi}_j^{cons}(\mathbf{p},\boldsymbol{\rho}_j) = \max_{\mathbf{y}_j \in cl[\tilde{Y}_j(\boldsymbol{\rho}_j)]} \pi_j^{cons}(\mathbf{p},\mathbf{y}_j), \quad \tilde{\pi}_i^{prod}(\mathbf{p},\boldsymbol{\varepsilon}_i) = \max_{\mathbf{x}_i \in cl[\tilde{X}_i(\boldsymbol{\varepsilon}_i)]} \pi_i^{prod}(\mathbf{p},\mathbf{x}_i). \quad (26)$$

(We note that $cl[\tilde{X}_i(\boldsymbol{\varepsilon}_i)] \subset X_i$, $cl[\tilde{Y}_j(\boldsymbol{\rho}_j)] \subset Y_j$.) Let $\tilde{P}^+(\boldsymbol{\varepsilon},\boldsymbol{\rho})$ denote the set of minimizers for (25). For a given $\tilde{\mathbf{p}}^+ \in \tilde{P}^+(\boldsymbol{\varepsilon},\boldsymbol{\rho})$, let $\tilde{\mathbf{x}}_i^+$, $\tilde{\mathbf{y}}_j^+$ denote the corresponding maximizers of the problems (26) and let us define $\tilde{\pi}_i^{*prod}(\boldsymbol{\varepsilon},\boldsymbol{\rho}) = \pi_i^{prod}(\tilde{\mathbf{p}}^+,\mathbf{x}_i^*)$, $\tilde{\pi}_i^{+prod}(\boldsymbol{\varepsilon},\boldsymbol{\rho}) = \pi_i^{prod}(\tilde{\mathbf{p}}^+,\tilde{\mathbf{x}}_i^+)$, $\tilde{\pi}_j^{*cons}(\boldsymbol{\varepsilon},\boldsymbol{\rho}) = \pi_j^{cons}(\tilde{\mathbf{p}}^+,\mathbf{y}_j^*)$, $\tilde{\pi}_j^{+cons}(\boldsymbol{\varepsilon},\boldsymbol{\rho}) = \pi_j^{cons}(\tilde{\mathbf{p}}^+,\tilde{\mathbf{y}}_j^+)$.

The duality gap is given by

$$\tilde{U}^D(\boldsymbol{\varepsilon},\boldsymbol{\rho}) - \tilde{U}(\boldsymbol{\varepsilon},\boldsymbol{\rho}) = \sum_{j \in J}[\tilde{\pi}_j^{+cons}(\boldsymbol{\varepsilon},\boldsymbol{\rho}) - \tilde{\pi}_j^{*cons}(\boldsymbol{\varepsilon},\boldsymbol{\rho})] + \sum_{i \in I}[\tilde{\pi}_i^{+prod}(\boldsymbol{\varepsilon},\boldsymbol{\rho}) - \tilde{\pi}_i^{*prod}(\boldsymbol{\varepsilon},\boldsymbol{\rho})].$$

Since $\tilde{X}_i(\boldsymbol{\varepsilon}_i) \subset X_i$, $\tilde{Y}_j(\boldsymbol{\rho}_j) \subset Y_j$, we have $\tilde{U}^D(\boldsymbol{\varepsilon},\boldsymbol{\rho}) \leq U^D$. Hence,

$$U = \tilde{U}(\boldsymbol{\varepsilon},\boldsymbol{\rho}) \leq \tilde{U}^D(\boldsymbol{\varepsilon},\boldsymbol{\rho}) \leq U^D,$$

which entails the relation between duality gaps of the original and the modified optimization problems:

$$0 \leq \tilde{U}^D(\boldsymbol{\varepsilon},\boldsymbol{\rho}) - \tilde{U}(\boldsymbol{\varepsilon},\boldsymbol{\rho}) \leq U^D - U, \forall \boldsymbol{\varepsilon} \in R_{>0}^{nT}, \forall \boldsymbol{\rho} \in R_{>0}^{mT}.$$

As a result, the total uplift needed to support the centralized dispatch solution $\forall \tilde{p}^+ \in \tilde{P}^+(\boldsymbol{\varepsilon},\boldsymbol{\rho})$ in the modified optimization problem (24) is no higher than that for the original problem (6), $\forall p^+ \in P^+$. We propose considering the limit as $\boldsymbol{\varepsilon} \to +0$, $\boldsymbol{\rho} \to +0$, and utilize $\tilde{P}^+(+0,+0)$ as the set of market prices and the corresponding lost profits $[\tilde{\pi}_j^{+cons}(+0,+0) - \tilde{\pi}_j^{*cons}(+0,+0)]$ for a consumer $j$ and $[\tilde{\pi}_i^{+prod}(+0,+0) - \tilde{\pi}_i^{*prod}(+0,+0)]$ for a producer $i$ as the individual uplift payments. Likewise, if according to the market procedures the consumer's payment for the volumes associated with (some of) the sunk benefit cannot exceed certain amount (for example, because of the activities initiated in the previous market planning period and the consumer's intertemporal constraints), then the corresponding compensation is not included in the considered uplift.

For the price-sensitive consumer bids, (20) is supplemented by the analogous definitions for the demand side, and we have the following statement similar to Proposition 1.

*Proposition 2:* If for some $j$ after optimization over the binary variables $\mathbf{v}_j$ the function $B_j$ becomes concave function of $\mathbf{d}_j$, then a price $\mathbf{p}$ supports solution $\mathbf{y}_j^*$ in the decentralized dispatch problem iff it supports solution $\mathbf{y}_j^*$ in the modified decentralized dispatch problem, $\forall \boldsymbol{\rho}_j \in R_{>0}^T$:

$$\mathbf{y}_j^* \in \arg\max_{\mathbf{y}_j \in Y_j}[B_j(\mathbf{y}_j) - \mathbf{p}^T\mathbf{d}_j] \Leftrightarrow \mathbf{y}_j^* \in \arg\max_{\mathbf{y}_j \in cl[\tilde{Y}_j(\boldsymbol{\rho}_j)]}[B_j(\mathbf{y}_j) - \mathbf{p}^T\mathbf{d}_j].$$

*Proof.* The consideration is fully analogous to the proof of Proposition 1. We also note that Proposition 2 stays valid if $\mathbf{y}_j^*$ is replaced by $\forall \overline{\mathbf{y}}_j \in \overline{\Omega}_j^{cons}$.

Thus, if no non-convexities are present in the power system, then both the convex hull pricing and the modified convex hull pricing algorithms result in the identical sets of prices (i.e. marginal prices) with zero total uplifts.

## 6. Applications and examples

In this section (except for Example 9), for simplicity, we restrict our consideration to the single-period power market. This model is attractive as it has exact analytic solutions in simple cases and allows clarifying the essential features of the proposed method. An $i$-th generator is assumed to have an output cost function $C_i(x_i) = c_i(g_i) + w_i u_i$ with variable cost function $c_i(g_i)$ and the fixed cost $w_i$, $w_i \geq 0$. The function $C_i(x_i)$ is defined on $X_i = \{x_i \mid u_i \in \{0,1\}, g_i \in R_{\geq 0}, u_i g_i^{\min} \leq g_i \leq u_i g_i^{\max}\}$, where $g_i^{\min}/g_i^{\max}$ denotes minimum/maximum output limits of the generator $i$. Also, $c_i(g_i)$ is assumed to be continuous non-decreasing convex continuous function on $\{0\} \cup [g_i^{\min}, g_i^{\max}]$ with $c_i(0) = 0$. Initially, all the generating units are assumed to be offline.

Let us define minimum economic output $g_i^{ec.\min}$ as follows: if $w_i = 0$, then $g_i^{ec.\min} = g_i^{\min}$; if $w_i \neq 0$, then $g_i^{ec.\min}$ is the lowest solution (if any) to $[w_i + c_i(g_i^{ec.\min})] \leq g_i^{ec.\min} \partial_+ c_i(g_i^{ec.\min})$ for $g_i^{\min} \leq g_i^{ec.\min} < g_i^{\max}$; if there is no solution, then $g_i^{ec.\min} = g_i^{\max}$. Hence, $g_i^{ec.\min}$ depends on the generator cost function and its individual feasible set $X_i$ only. It is straightforward to obtain the following expression for the subdifferential $\partial \pi_i^{prod}(p)$:

$$\partial \pi_i^{prod}(p) = \begin{cases} 0, p < [w_i + c_i(g_i^{ec.\min})]/g_i^{ec.\min} \\ [0, \gamma_{i\max}], p = [w_i + c_i(g_i^{ec.\min})]/g_i^{ec.\min} \\ \Gamma_i(p), p > [w_i + c_i(g_i^{ec.\min})]/g_i^{ec.\min} \end{cases}, (27)$$

where $\gamma_{i\max}$ is the maximum output volume $g_i$ satisfying $[w_i + c_i(g_i)] \leq g_i \partial_+ c_i(g_i)$ for $g_i^{\min} \leq g_i \leq g_i^{\max}$, and $\Gamma_i(p)$ denotes elements (which may depend on $p$) no lower than $\gamma_{i\max}$. We note that if $g_i^{ec.\min} = 0$, i.e. both $w_i = 0$ and $g_i^{\min} = 0$, then (27) is well-defined because $\lim_{g_i \to +0} c_i(g_i)/g_i = \partial_+ c_i(0)$, which is finite.

The set of maximizers of (2) has the output volumes, which we denote as $G_i(p)$, given by

$$G_i(p) = \begin{cases} 0, p < [w_i + c_i(g_i^{ec.\min})]/g_i^{ec.\min} \\ \{0\} \cup \{\gamma_{i\max}\}, p = [w_i + c_i(g_i^{ec.\min})]/g_i^{ec.\min} \\ \Gamma_i(p), p > [w_i + c_i(g_i^{ec.\min})]/g_i^{ec.\min} \end{cases}. (28)$$

Clearly, $G_i(p)$ is the individual supply curve of the unit $i$. Note that $G_i(p) \subset \{0\} \cup [g_i^{ec.\min}, g_i^{\max}]$, $\forall p$. Thus, if $g_i^{ec.\min} \neq 0$, i.e. if $w_i \neq 0$ and/or $g_i^{\min} \neq 0$, then the supply curve has a gap: under no value of $p$ the output volumes from the range $(0, g_i^{ec.\min})$ are supplied in the decentralized dispatch problem. That is the essence of the difference between the optimal generator outputs in the primal and the decentralized dispatch problems (and also the problem dual to (1)) because, for a unit with $g_i^{\min} < g_i^{ec.\min}$, the optimal generator output in the primal problem can belong to $(g_i^{\min}, g_i^{ec.\min})$, while the optimal generator output in the decentralized dispatch problem cannot belong to that range. From (27) and (28) we deduce $\partial_{\pm}\pi_i^{prod}(p) \in G_i(p)$, and, therefore, the output volumes $\partial_{\pm}\pi_i^{prod}(p)$ at a price $p$ belong to the supply curve of the unit $i$.

The binary variable can be excluded at the expense of having discontinuity introduced in the cost function: let $f_i(g_i) = w_i\theta(g_i) + c_i(g_i)$ with $dom[f_i(g_i)] = \{0\} \cup [g_i^{min}, g_i^{max}]$ and step-function $\theta(g_i)$ defined equal to 1 for $x_i > 0$ and 0 otherwise. Let us further define $f_i(g_i) = +\infty$ for $g_i \notin dom[f_i(g_i)]$ and denote by $f_i^h(g_i)$ the closed convex hull of the extended-value function $f_i(g_i)$. Thus, $\pi_i^{prod}(p)$ is Fenchel conjugate of $f_i(g_i)$, which implies [29]

$$\pi_i^{prod}(p) = \max_{g_i \in dom[f_i^h(g_i)]} [pg_i - f_i^h(g_i)]. \quad (29)$$

We note that in the general case the set of output volumes that maximize (2) is different from the maximizers of (29). Also, the inversion rule for the subgradient relations [30] yields

$$g_i \in \partial\pi_i^{prod}(p) \leftrightarrow p \in \partial f_i^h(g_i).$$

Thus, the set of points $(g_i, \partial f_i^h(g_i))$ would be the supply curve for the unit $i$ if its cost function were $f_i^h(g_i)$. This observation simplifies analysis in the numerous examples below.

## 6.1 The uninode single-period power systems with fixed load (no price-sensitive consumer bids)

Let us consider the case of the uninode single-period power system with fixed demand $d$ and zero minimum output limits for all units: $g_i^{min} = 0$, $\forall i \in I$. Such systems were extensively studied in [28], below we give brief summary of that analysis. The set $\overline{\Omega}$ is given by the feasible set of the primal problem excluding the points with at least one generating unit with nonzero start-up cost having zero output volume and a status "ON". The set $\overline{\Omega}_i^{prod}$, $\forall i \in I$, consists of elements with volumes $g_i \in [\overline{g}_i^{min}, \overline{g}_i^{max}]$, where $\overline{g}_i^{min} = \max(d - \sum_{i':i' \in I, i' \neq i} g_{i'}^{max}; 0)$, $\overline{g}_i^{max} = \min(d; g_i^{max})$, and the corresponding values of the status variable $u_i$. We note that if both $\overline{g}_i^{min} = 0$ and $w_i \neq 0$ for some unit $i$, then the set $\overline{\Omega}_i^{prod}$ is not compact (when $d > 0$). We have the following explicit expressions for the sets $\bigcup_{\overline{x}_i \in \overline{\Omega}_i^{prod}} \Delta_i^{prod}(\overline{x}_i, \varepsilon_i)$, $\forall i \in I$, $\forall \varepsilon_i > 0$, [28]:

- if $\overline{g}_i^{min} > \varepsilon_i$, then
  $\bigcup_{\overline{x}_i \in \overline{\Omega}_i^{prod}} \Delta_i^{prod}(\overline{x}_i, \varepsilon_i) = \{(u_i, g_i) \mid u_i = 1, \overline{g}_i^{min} - \varepsilon_i \leq g_i \leq \min[\overline{g}_i^{max} + \varepsilon_i; g_i^{max}]\}$;
- if $\overline{g}_i^{min} \leq \varepsilon_i$, then
  $\bigcup_{\overline{x}_i \in \overline{\Omega}_i^{prod}} \Delta_i^{prod}(\overline{x}_i, \varepsilon_i) = \{(u_i, g_i) \mid u_i = 1, 0 \leq g_i \leq \min[\overline{g}_i^{max} + \varepsilon_i; g_i^{max}]\}$.

We note that to guarantee the non-confiscatory pricing for power it is sufficient to ensure that $\{(0,0)\} \in \widetilde{X}_i(\varepsilon_i)$ even if $\Psi_i^{prod}$ contains some other elements as well, [28]. However, to illustrate application of the proposed pricing method, below we adopt the general expression (12). We also note that the resulting sets $\widetilde{X}_i(\varepsilon_i)$ are compact, $\forall \varepsilon_i > 0$, $\forall i \in I$. It was shown in [28] that $\overline{g}_i^{\min}$, $\forall i \in I$, are irrelevant for both the price and the uplift calculations, and they can be formally set to zero for all the generating units. Let us identify the generating units with $g_i^{ec.\min}$ exceeding demand (named LNMGUs in [28]). For all units, except for LNMGUs, the values of $\overline{g}_i^{\max}$ are irrelevant as well, and the original sets $X_i$ can be used instead of $\widetilde{X}_i(\varepsilon_i)$ in the modified convex hull pricing method with no effect on the resulting set of the market prices and the uplift payment for each generating unit. Hence, only LNMGUs require special treatment in the modified approach compared to the convex hull pricing mechanism.

Now, let us turn to the power systems with fixed load and possibly nonzero minimum output limits. In this case, the proposed pricing scheme may result in a different pricing outcome compared to the convex hull pricing method even for the power systems without LNMGUs, as shown in the example below.

*Example 3.*

Consider a power system with fixed demand $d=200$ MWh and two generating units having $C_i(x_i) = a_i g_i + w_i u_i$ with parameter values specified in Table 1.

Table 1. Parameters values for Example 3.

|  | $g_i^{\min}$, MWh | $g_i^{\max}$, MWh | $a_i$, \$/MWh | $w_i$, \$ |
|---|---|---|---|---|
| Unit 1 | 80 | 160 | 20.00 | 0.00 |
| Unit 2 | 80 | 160 | 30.00 | 15.00 |

Clearly, the primal problem has the unique solution $g_1^* = 120 MWh$, $g_2^* = 80 MWh$, $u_1^* = u_2^* = 1$. The convex hull pricing results in the unique market price $p^+ = a_2 + w_2/g_2^{\max} = \$30.09/MWh$ and $g_1^+ = g_2^+ = 160 MWh$, $u_1^+ = u_2^+ = 1$, (we note that there is also the other solution for the dual problem: $g_1^+ = 160 MWh$, $u_1^+ = 1$, $u_2^+ = g_2^+ = 0$). The total uplift in the convex hull pricing method is given by \$411.40.

Now, let us apply the modified convex hull pricing framework. Clearly, both generators have to be online to satisfy the fixed demand, hence $\overline{\Omega}_i^{prod} = \{(u_i, g_i) | u_i = 1, g_i \in [\overline{g}_i^{\min}, \overline{g}_i^{\max}])\}$ with $\overline{g}_i^{\min} = g_i^{\min} = 80 MWh$, $i=1,2$, $\overline{g}_1^{\max} = d - g_2^{\min} = 120 MWh$, $\overline{g}_2^{\max} = d - g_1^{\min} = 120 MWh$. For sufficiently small $\varepsilon_i > 0$, we have $\widetilde{X}_i(\varepsilon_i) = \Psi_i^{prod} \cup \{(u_i, g_i) | u_i = 1, g_i \in [\overline{g}_i^{\min}, \overline{g}_i^{\max} + \varepsilon_i]\}$, $i=1,2$, with $\Psi_1^{prod} = \{(0,0)\} \cup \{(1,0)\}$, $\Psi_2^{prod} = \{(0,0)\}$, and the unique market price $\widetilde{p}^+(\varepsilon_1, \varepsilon_2) = a_2 + w_2/(\overline{g}_2^{\max} + \varepsilon_2)$, which gives $\widetilde{p}^+(+0,+0) = \$30.13/MWh$ and the total uplift payment equals \$4.60, payable to generator 2. Details of the comparison of two methods are given in Table 2.

Table 2. The uplift payments for Example 3.

|  | Convex hull pricing | | | Modified convex hull pricing | | |
|---|---|---|---|---|---|---|
|  | $\pi_i^{*prod}$, \$ | $\pi_i^{+prod}$, \$ | Uplift, \$ | $\widetilde{\pi}_i^{*prod}$, \$ | $\widetilde{\pi}_i^{+prod}$, \$ | Uplift, \$ |
| Unit 1 | 1210.80 | 1614.40 | 403.60 | 1215.60 | 1215.60 | 0.00 |
| Unit 2 | -7.80 | 0.00 | 7.80 | -4.60 | 0.00 | 4.60 |
| Total | 1203.00 | 1614.40 | 411.40 | 1211.00 | 1215.60 | 4.60 |

Another specific property of the power systems with the fixed load and nonzero minimum output limits is that the set $\widetilde{X}_i(\varepsilon_i)$ generally cannot be replaced by $X_i$ even for the units with zero minimum output limits and start-up costs - this is illustrated in the following example.

*Example 4.*

Let us amend the power system described in Example 3 by setting $g_1^{\min} = 0 MWh$. This substitution changes neither the primal problem solution nor the convex hull pricing outcome. Application of the modified convex hull pricing framework yields $\overline{\Omega}_i^{prod} = \{(u_i, g_i) | u_i = 1, g_i \in [\overline{g}_i^{\min}, \overline{g}_i^{\max}])\}$, $i = 1, 2$, with $\overline{g}_1^{\min} = d - g_2^{\max} = 40 MWh$, $\overline{g}_1^{\max} = d - g_2^{\min} = 120 MWh$, $\overline{g}_2^{\min} = g_2^{\min} = 80 MWh$, $\overline{g}_2^{\max} = 160 MWh$. For sufficiently small $\varepsilon_1$ we have

$$\widetilde{X}_1(\varepsilon_1) = \{(0,0)\} \cup \{(1,0)\} \cup \{(u_1, g_1) | u_1 = 1, g_1 \in [\overline{g}_1^{\min} - \varepsilon_1, \overline{g}_1^{\max} + \varepsilon_1]\},$$
$$\widetilde{X}_2(\varepsilon_2) = \{(0,0)\} \cup \{(u_2, g_2) | u_2 = 1, g_2 \in [\overline{g}_2^{\min}, \overline{g}_2^{\max} + \varepsilon_2]\}.$$

The market price is unique and equal to $\widetilde{p}^+(\varepsilon_1, \varepsilon_2) = a_2 + w_2/(\overline{g}_2^{\max} + \varepsilon_2)$, which implies $\widetilde{p}^+(+0, +0) = \$30.09/MWh$ and the total uplift of \$7.80. In Table 3 we give the unit profits and uplifts in the modified convex hull pricing method.

**Table 3.** The uplift payments for Example 4 in the proposed method.

|  | Modified convex hull pricing | | |
|---|---|---|---|
|  | $\widetilde{\pi}_i^{*prod}$,\$ | $\widetilde{\pi}_i^{+prod}$,\$ | Uplift,\$ |
| Unit 1 | 1210.80 | 1210.80 | 0 |
| Unit 2 | -7.80 | 0.00 | 7.80 |
| Total | 1203.00 | 1210.80 | 7.80 |

Replacement of $\widetilde{X}_i(\varepsilon_i)$ by $X_1$ in the proposed pricing algorithm produces the same price $\widetilde{p}^+(+0,+0)$ and uplift payment for the second unit but different uplift payment for the first unit, namely, $(\widetilde{p}^+ - a_1)(g_1^{\max} - g_1^*) = \$1203.60$. Thus, although the first generating unit has the convex cost function (both zero minimum output limit and zero start-up cost as well as the convex variable cost), utilization of the unit's original individual feasible set $X_1$ in the proposed pricing method results in the different outcome.

**6.2 The uninode single-period power systems with price-sensitive consumer bids**

For simplicity, we consider the uninode single-period power system with consumers having the concave benefit functions $B_j(d_j)$ of the power consumption volumes $d_j$, $d_j \in Y_j$, $Y_j = \{d_j | d_j \in R_{\geq 0}, d_j^{\min} \leq d_j \leq d_j^{\max}\}$, $\forall j \in J$, with constant parameters $d_j^{\min}$ and $d_j^{\max}$ satisfying $0 \leq d_j^{\min} \leq d_j^{\max} < +\infty$ and representing fixed part of a consumer $j$ demand volume and its maximal consumption volume, respectively.

First, we consider the case of all generating units in the power system having zero minimum output limits. To find explicit expression for $\overline{\Omega}_i^{prod}$, we consider the maximum value of a generator $i$ economically and technologically feasible output, which we denote as $\overline{g}_i^{\max}$.

*Proposition 3*: If $g_i^{\min}=0$, $\forall i \in I$, then $\forall i_0 \in I$:

- if $\sum_{j \in J} d_j^{\min} < g_{i_0}^{\max}$, then $\bar{g}_{i_0}^{\max}$ is the optimal generator $i_0$ output volume determined by solving the following problem (if the problem has multiple optimal points, then $\bar{g}_{i_0}^{\max}$ corresponds to the one with the maximum generator $i_0$ output):

$$\max_{\substack{x_{i_0}, d, \\ x_{i_0} \in X_{i0}, \\ d_j \in Y_j, \forall j \in J \\ g_{i_0} = \sum_{j \in J} d_j}} \left[ \sum_{j \in J} B_j(d_j) - C_{i_0}(x_{i_0}) \right].$$

- if $\sum_{j \in J} d_j^{\min} \geq g_{i_0}^{\max}$, then $\bar{g}_{i_0}^{\max} = g_{i_0}^{\max}$.

*Proof*. The specific feature of the uninode single-period system is that each generator competes with the other generators and doesn't need to cooperate with them to achieve its maximum output volume (except for the case when the fixed load segment of demand needs to be supplied): if operating statuses "ON"/"OFF" of all units are fixed, then a decrease in some unit output will either increase or leave unchanged the output of any other unit. Thus, in the case of $\sum_{j \in J} d_j^{\min} < g_{i_0}^{\max}$, we may start with the collection of $g_i^{\max}$, $i \in I$, reduce some of $g_i^{\max}$ (except for the unit $i_0$), provided that the resulting centralized dispatch problem is feasible, and compute the centralized dispatch outcome with the new feasible set. Let us define the algorithm: for all units (except for the unit $i_0$) set $g_i^{\max}$ equal to the new centralized dispatch output volumes, then reduce some of $g_i^{\max}$, except for the unit $i_0$, provided that the resulting centralized dispatch is feasible, find the solution to the optimization problem with the new limits on the maximum output. If the algorithm is run repeatedly, then at each step the output of unit $i_0$ either increases or stays the same. Thus, the maximum possible output of the generator $i_0$ is realized when no other generator is operating, which results in the statement of the first bullet. If $\sum_{j \in J} d_j^{\min} \geq g_{i_0}^{\max}$, i.e. the fixed load part of demand exceeds the generator $i_0$ maximum output, then we let all the other units produce exactly $(\sum_{j \in J} d_j^{\min} - g_{i_0}^{\max})$. To achieve that, the generator $i_0$ may sign the bilateral contracts for power with each generating unit $i$, $i \neq i_0$, to supply $\bar{g}_i = g_i^* (\sum_{j \in J} d_j^{\min} - g_{i_0}^{\max}) / \sum_{i':i' \in I, i' \neq i_0} g_{i'}^*$, engage in the contracts with the fixed load segment of demand to resell these power volumes, and sign the contracts with the remaining part of the fixed load to supply $g_{i_0}^{\max}$. It is straightforward to check that $0 \leq \bar{g}_i \leq g_i^{\max}$. Therefore, the set $\{\bar{g}_i\}$ is realized as the centralized dispatch outcome with only fixed load segment of demand cleared. This gives the statement of the second bullet. The proposition is proved.

*Proposition 4*: If $g_i^{\min} = 0$, $\forall i \in I$, then $\forall i_0 \in I$:

- if $\sum_{j \in J} d_j^{\min} = 0$, $\bar{g}_{i_0}^{\max} > 0$, and $w_{i_0} = 0$, then $\bar{\Omega}_{i_0}^{prod} = \{(0,0)\} \cup \{(u_{i_0}, g_{i_0}) \mid u_{i_0} = 1, g_{i_0} \in [0, \bar{g}_{i_0}^{\max}]\}$;

- if $\sum_{j \in J} d_j^{\min} = 0$, $\bar{g}_{i_0}^{\max} > 0$, and $w_{i_0} \neq 0$, then $\overline{\Omega}_{i_0}^{prod} = \{(0,0)\} \cup \{(u_{i_0}, g_{i_0}) \mid u_{i_0} = 1, g_{i_0} \in [\bar{g}_{i_0}^{\min}, \bar{g}_{i_0}^{\max}]\}$, where $\bar{g}_{i_0}^{\min}$ denotes the minimum output, satisfying

$$c_{i_0}(\bar{g}_{i_0}^{\min}) + w_{i_0} = \max_{\substack{d, \\ d_j \in Y_j, \forall j \in J \\ \bar{g}_{i_0}^{\min} = \sum_{j \in J} d_j}} \sum_{j \in J} B_j(d_j); \qquad (30)$$

- if $\sum_{j \in J} d_j^{\min} = 0$, $\bar{g}_{i_0}^{\max} = 0$, and $w_{i_0} = 0$, then $\overline{\Omega}_{i_0}^{prod} = \{(0,0)\} \cup \{(1,0)\}$;
- if $\sum_{j \in J} d_j^{\min} = 0$, $\bar{g}_{i_0}^{\max} = 0$, and $w_{i_0} \neq 0$, then $\overline{\Omega}_{i_0}^{prod} = \{(0,0)\}$;
- if $\sum_{j \in J} d_j^{\min} > 0$, then $\overline{\Omega}_{i_0}^{prod}$ is comprised of all elements with the output volume $g_{i_0} \in [\bar{g}_{i_0}^{\min}, \bar{g}_{i_0}^{\max}]$, where $\bar{g}_{i_0}^{\min} = \max(\sum_{j \in J} d_j^{\min} - \sum_{i: i \in I, i \neq i_0} g_i^{\max}; 0)$, and the corresponding statuses (excluding the element $\{(1,0)\}$, if any, in case of $w_{i_0} \neq 0$).

*Proof.* The proof is straightforward and is based on the observation that if $\sum_{j \in J} d_j^{\min} = 0$, then we may start with the unit $i_0$ output $\bar{g}_{i_0}^{\max}$ and, keeping all the other units offline, gradually reduce $\bar{g}_{i_0}^{\max}$. If $\sum_{j \in J} d_j^{\min} > 0$, then we also start with the output $\bar{g}_{i_0}^{\max}$ and gradually reduce $\bar{g}_{i_0}^{\max}$ with the other units utilized only to ensure that part of (or all) the fixed load segment of demand is satisfied by their output volumes. The proposition is proved.

Contrary to competing producers, a consumer may find it mutually beneficial to cooperate with the other consumers to increase its consumption volume: if the start-up cost is present, reduction of the consumption volume by one consumer may lower the consumption volume of another consumer resulting from the centralized dispatch. Hence, the consumer may find it profitable to conclude bilateral contracts with some of the other consumers even with higher bids. Thus, Propositions 3 and 4 are not trivially generalized for $\overline{\Omega}_j^{cons}$. We also note that $\overline{\Omega}_i^{prod}$ and $\overline{\Omega}_j^{cons}$ may not be compact.

Now we illustrate application of the modified convex hull pricing method for Examples 1 and 2 and show that the total uplift needed to provide economic stability of the centralized dispatch outcome is reduced compared to the total uplift implied by the convex hull pricing algorithm.

*Example 1 (revisited).*
Since no nonzero power volume can be fully paid for by the consumer, we have $\overline{\Omega}^{prod} = \{(0,0)\}$, $\overline{\Omega}^{cons} = \{0\}$. Therefore,

$$\widetilde{X}(\varepsilon) = \{(0,0)\} \cup \{(u,g) \mid u = 1, g \in [0, \varepsilon]\}, \ \widetilde{Y}(\rho) = [0, \rho].$$

For $0 < \varepsilon < w/(b-a)$, we have $\widetilde{P}^+(\varepsilon, \rho) = [b, a + w/\varepsilon]$, and

$$\widetilde{\pi}^{+prod}(\varepsilon) = \widetilde{\pi}^{+cons}(\rho) = \pi^{*prod}(\varepsilon) = \pi^{*cons}(\rho) = 0, \ \forall \widetilde{p}^+ \in \widetilde{P}^+(\varepsilon, \rho).$$

Hence, no market player receives any compensation and the total uplift payment is zero. In the limit as $\varepsilon \to +0$, $\rho \to +0$, we have $\widetilde{P}^+(+0, +0) = [b, +\infty)$.

*Example 2 (revisited).*
Let us apply the proposed approach to Example 2 under assumption that there is only one consumer in the power system. Application of Proposition 4 gives $\overline{\Omega}^{prod} = \{(0,0)\} \cup \{(u,g) \mid u = 1, g \in [\bar{g}^{\min}, d^*]\}$ with $\bar{g}^{\min} = d^{\max}[1 - \sqrt{1 - 4w/(ad^{\max})}]/2$.

The power balance constraint gives $\overline{\Omega}^{cons} = \{0\} \cup [\overline{g}^{\min}, d^*]$. Hence, for sufficiently small $\varepsilon > 0$ and $\rho > 0$, we have $\widetilde{X}(\varepsilon) = \{(0,0)\} \cup \{(u,g) | u = 1, g \in [\overline{g}^{\min} - \varepsilon, d^* + \varepsilon]\}$ for the consumer and $\widetilde{Y}(\rho) = \{0\} \cup [\overline{g}^{\min} - \rho, d^* + \rho]$ for the producer. It is straightforward to find that the market price is unique and given by $\widetilde{p}^+(\varepsilon, \rho) = a + w/(g^* + \varepsilon)$, which entails that in the limit as $\varepsilon \to +0$, $\rho \to +0$ we have $\widetilde{\pi}^{+prod} = \widetilde{\pi}^{*prod} = 0$. Consequently, no uplift is paid to the generator, while the consumer's uplift equals $\widetilde{\pi}^{+cons} - \widetilde{\pi}^{*cons} = w^2/(ad^{\max})$. Therefore, $\widetilde{U}^D(0,0) - \widetilde{U}(0,0) = w^2/(ad^{\max})$. Using the stated assumptions on the parameters, it is straightforward to check that $\widetilde{U}^D(0,0) - \widetilde{U}(0,0) \leq w/6$, while $U^D - U > w\left(1/2 + wd^{\max}/[(2g^{\max})^2 a]\right) > w/2$. Therefore, the total uplift produced by the proposed method is lower than that implied by the convex hull pricing mechanism.

The convex hull pricing approach produces the outcome that is independent of the structure of demand: substitution of any group of consumers by one consumer with the aggregate bid will not produce a different set of the market prices or a total uplift. The proposed modified convex hull pricing method lacks this property as the total uplift may depend on the structure of demand as shown in the example below, which also illustrates the application of the proposed method for power systems with nonzero minimum output limits and price-sensitive load.

*Example 5.*

Consider a power system with two consumers having the benefit functions $B_j(d_j) = b_j d_j$, $d_j \in Y_j$, $Y_j = [0, d_j^{\max}]$, $j = 1, 2$, with parameter values specified in Table 4.

**Table 4.** Parameter values for Example 5.

|  | $b_j, \$/MWh$ | $d_j^{\max}, MWh$ |
|---|---|---|
| Consumer 1 | 100.00 | 100 |
| Consumer 2 | 15.00 | 300 |

Let us introduce one generating unit with $g^{\min} = g^{\max} = 250 MWh$ and cost function $C(x) = ag + wu$ with $a = \$20.00/MWh$ and $w = \$50.00$. Clearly, $X = \{(0,0)\} \cup \{(1, 250)\}$. It is straightforward to check that the model parameters satisfy $d_1^{\max} < g^{\max} < d_2^{\max}$, $ag^{\max} + w < b_1 d_1^{\max} + b_2(g^{\max} - d_1^{\max})$. Thus, the generator is either online operating at maximum output or offline. The primal problem has the unique solution given by $d_1^* = d_1^{\max} = 100 MWh$, $d_2^* = g^{\max} - d_1^{\max} = 150 MWh$, $u^* = 1$, $g^* = 250 MWh$. We have $\overline{\Omega}^{prod} = X$ for the generator, while for consumers these sets are $\overline{\Omega}_1^{cons} = \{0\} \cup [\overline{d}_1^{\min}, d_1^{\max}]$ and $\overline{\Omega}_2^{cons} = \{0\} \cup [g^{\max} - d_1^{\max}, g^{\max} - \overline{d}_1^{\min}]$ with $\overline{d}_1^{\min} = [(a - b_2)g^{\max} + w]/(b_1 - b_2)$. For any $\varepsilon > 0$ and sufficiently small positive $\rho_1$ and $\rho_2$, we have

$$\widetilde{X}(\varepsilon) = X, \widetilde{Y}_1(\rho_1) = [0, \rho_1] \cup [\overline{d}_1^{\min} - \rho_1, d_1^{\max}],$$
$$\widetilde{Y}_2(\rho_2) = [0, \rho_2] \cup [g^{\max} - d_1^{\max} - \rho_2, g^{\max} - \overline{d}_1^{\min} + \rho_2].$$

The modified convex hull pricing procedure yields the unique price $\tilde{p}^+(+0,+0,+0) = a + w/g^{max} = \$20.20/MWh$ and the profits and individual uplifts presented in Table 5.

Table 5. The uplift payments for Example 5 in the proposed method.

|  | Modified convex hull pricing | | |
|---|---|---|---|
|  | $\tilde{\pi}^*$,$ | $\tilde{\pi}^+$,$ | Uplift,$ |
| Consumer 1 | 7980.00 | 7980.00 | 0.00 |
| Consumer 2 | -780.00 | 0.00 | 780.00 |
| Producer | 0.00 | 0.00 | 0.00 |
| Total | 7200.00 | 7980.00 | 780.00 |

Therefore, both the generator and consumer 1 uplifts are zero, while consumer 2 uplift is $780.00, which is the total uplift in this case. We note that application of the convex hull pricing algorithm results in identical market price and individual uplifts.

Now, let us assume that the same aggregate demand is represented by just one consumer, then its benefit function is $B(d) = b_1 \min(d; d_1^{max}) + b_2 \max(d - d_1^{max}; 0)$, $0 \le d \le d_1^{max} + d_2^{max}$. In this case, we have $\overline{\Omega}^{cons} = \{0\} \cup \{g^{max}\}$ with the same as above expressions for $\overline{\Omega}^{prod}$, $\tilde{X}(\varepsilon)$, and $\tilde{p}^+$. For sufficiently small positive $\rho$, we have $\tilde{Y}(\rho) = [0, \rho] \cup [g^{max} - \rho, g^{max} + \rho]$, which entails

$$\tilde{\pi}^{+cons}(\varepsilon, \rho) = b_1 d_1^{max} + b_2(g^{max} - \rho - d_1^{max}) - \tilde{p}^+(g^{max} - \rho),$$
$$\tilde{\pi}^{*cons}(\varepsilon, \rho) = b_1 d_1^{max} + b_2(d_1^{max} - g^{max}) - \tilde{p}^+ g^{max}.$$

In the limit as $\varepsilon \to +0$, $\rho \to +0$ we have the pricing outcomes summarized in Table 6.

Table 6. The uplift payments for Example 5 in the case of one consumer.

|  | Convex hull pricing | | | Modified convex hull pricing | | |
|---|---|---|---|---|---|---|
|  | $\pi^*$,$ | $\pi^+$,$ | Uplift,$ | $\tilde{\pi}^*$,$ | $\tilde{\pi}^+$,$ | Uplift,$ |
| Consumer | 7200.00 | 7980.00 | 780.00 | 7200.00 | 7200.00 | 0.00 |
| Producer | 0.00 | 0.00 | 0.00 | 0.00 | 0.00 | 0.00 |
| Total | 7200.00 | 7980.00 | 780.00 | 7200.00 | 7200.00 | 0.00 |

Therefore, the total uplift payment is zero in the modified convex hull pricing approach. Thus, we conclude that the total uplift payment in the modified convex pricing approach depends on the structure of demand: knowledge of the aggregate demand is insufficient for the uplift calculation and the individual consumer benefit functions are needed.

We also note the following property of the proposed method. If some generating unit $i_0$ with $w_{i_0} \ne 0$ and/or $g_{i_0}^{min} \ne 0$ has both $\overline{\Omega}_{i_0}^{prod}$ and $\Psi_{i_0}^{prod}$ comprised of elements with zero output volumes only, then in the limit as $\varepsilon_{i_0} \to +0$ this unit does not contribute to the set of market prices $\tilde{P}^+$ and its individual uplift vanishes. Indeed, if $g_{i_0}^{min} \ne 0$ and $w_{i_0} \ne 0$, then for sufficiently small positive $\varepsilon_{i_0}$ we have $\tilde{X}_{i_0}(\varepsilon_{i_0}) = \{(0,0)\}$, and the unit's profit in the decentralized dispatch optimization problem is zero at any price. Likewise, if $g_{i_0}^{min} \ne 0$ and $w_{i_0} = 0$, then $\tilde{X}_{i_0}(\varepsilon_{i_0}) = \{(0,0)\} \cup \{(1,0)\}$ for sufficiently small positive $\varepsilon_{i_0}$ and the unit's profit in the decentralized dispatch optimization problem vanishes at any price. If $g_{i_0}^{\min} = 0$, but $w_{i_0} \neq 0$, then the unit contributes to objective function in dual of the modified primal problem only for prices above $[w_{i_0} + c_{i_0}(\varepsilon_{i_0})]/\varepsilon_{i_0} = w_{i_0}/\varepsilon_{i_0} + O(1)$, which goes to infinity as $\varepsilon_{i_0} \to +0$. Therefore, in case of uninode single-period market, such units can be removed from consideration in dual of the modified primal problem affecting neither a set the of market prices nor the uplifts in the limit as $\varepsilon_{i_0} \to +0$. In particular, this observation is relevant in the following extreme cases: if $g_{i_0}^{\min}$ of some generating unit exceeds the maximum consumption volume $\sum_{j \in J} d_j^{\max}$, likewise for a unit with the start-up cost higher than the total consumer benefit for the maximum consumption volume, - such units can be omitted from the consideration in the modified convex hull pricing algorithm for a uninode single-period market.

In case of no demand-side non-convexities, there is a question if the original consumer individual feasible sets (and, hence, the original demand curves) can be used to find the market price and the uplift payments in the proposed method. We will address this issue under the assumption of zero minimum output limits for all generating units. Let us construct the aggregate supply curve of the modified primal problem (24), which we refer to as the aggregate modified supply curve, and the aggregate demand curve of the primal problem (5). Each of the curves has at least one point for any value of the price. However, the former curve generally has gaps, while the latter is continuous. (We note, however, that the aggregate modified demand curve for (24) may also have gaps.) The following proposition entails that if the aggregate modified supply curve has a point with a price $p$ and a volume $q$, $q \geq \sum_{j \in J} d_j^{\min}$, such that the aggregate demand curve has a point with the volume $q$ and a price $p'$, $p' \geq p$, then the corresponding individual outputs (consumption) volumes are both economically and technologically feasible, in particular, the consumption volumes satisfy (14) with appropriate $\bar{x}$.

*Proposition 5*: Let $\hat{x} = (\hat{x}_1,..,\hat{x}_n)$, $\hat{x}_i = (\hat{u}_i, \hat{g}_i)$, $\hat{d} = (\hat{d}_1,..,\hat{d}_m)$ with

$$\hat{x}_i \in \arg \max_{x_i \in cl[\tilde{X}_i(\varepsilon_i)]} \pi_i^{prod}(p, x_i), \ \forall i \in I, \ \hat{d}_j \in \arg \max_{d_j \in Y_j} \pi_j^{cons}(p', d_j), \ \forall j \in J,$$

for some $p$, $p'$. If $g_i^{\min} = 0$, $\varepsilon_i > 0$, $\forall i \in I$, $p' \geq p$, and $\sum_{i \in I} \hat{g}_i = \sum_{j \in J} \hat{d}_j$, then $\{\hat{x}, \hat{d}\}$ satisfies (14).

*Proof.* Let $q = \sum_{i \in I} \hat{g}_i = \sum_{j \in J} \hat{d}_j$. Since the aggregate modified supply curve and the aggregate demand curve have points with equal volume at a price $p$ and $p'$ respectively with $p' \geq p$, the outputs $\hat{g}_i$ satisfy certain restrictions. Namely, if $\hat{g}_i < g_i^{ec.\min}$ with $g_i^{ec.\min}$ calculated using original generator's set $X_i$, then $\bar{g}_i^{\max} + \varepsilon_i < g_i^{ec.\min} \leq g_i^{\max}$ and, as a result, $\hat{g}_i = 0$ or $\hat{g}_i = \bar{g}_i^{\max} + \varepsilon_i$. The former case is trivial, while for the latter Proposition 3 implies that for $\varepsilon_i > 0$ all points on the aggregate demand curve with volume $\hat{g}_i$ have prices below the lowest price on the modified supply curve of the unit $i$ with volume $\hat{g}_i$. This contradicts the assumption of $p' \geq p$ given continuity and (weak) monotonicity of the aggregate demand curve. Therefore, $\forall i$ with $\hat{g}_i \neq 0$, we have $g_i^{ec.\min} \leq \hat{g}_i$. We also note that Proposition 3 yields $\hat{g}_i \leq \overline{g}_i^{\max}$. Hence, we conclude that if $\hat{g}_i \neq 0$, then $g_i^{ec.\min} \leq \hat{g}_i \leq \overline{g}_i^{\max}$. Consider the following optimization problem:

$$\max_{\substack{x,d, \\ x_i \in X_i, \forall i \in I \\ d_j \in Y_j, \forall j \in J \\ x \leq \hat{x}, d \leq \hat{d}}} \left[ \sum_{j \in J} B_j(d_j) - \sum_{i \in I} C_i(x_i) \right] \quad (31).$$

Now, we show that the conditions of this proposition imply that the aggregate supply curve for (31) has a point with the price $p$ and the volume $q$. To prove that, we demonstrate that for any producer $i$ the supply curve for (31) has a point with the price $p$ and the volume $\hat{g}_i$:

$$\hat{x}_i \in \arg \max_{\substack{x_i, \\ x_i \in X_i, \\ x_i \leq \hat{x}_i}} \pi_i^{prod}(p, x_i), \quad \forall i \in I.$$

Since $\hat{g}_i$ can be formally regarded as new value of $g_i^{\max}$, Proposition 4 classifies all the possible cases and its application to (31) yields the following.

If conditions of the 1-st bullet of Proposition 4 are met, then $g_i^{ec.\min} = 0$, $\hat{g}_i \leq \overline{g}_i^{\max}$, and we conclude that all points on the modified supply curve of the unit with volumes not exceeding $\hat{g}_i$ belong to the unit's supply curve for (31). In particular, it has a point with the price $p$ and the volume $\hat{g}_i$.

If conditions of the 2-nd bullet of Proposition 4 hold, then the following cases are possible: $\hat{g}_i = 0$ or $\overline{g}_i^{\min} - \varepsilon_i = \hat{g}_i$ or $\overline{g}_i^{\min} - \varepsilon_i < \hat{g}_i$. In the first case, the unit's supply curve for (31) has zero volumes for any price value and, hence, at the price $p$. The second case is not possible since $\overline{g}_i^{\min} - \varepsilon_i = \hat{g}_i$ entails $g_i^{ec.\min} \leq \overline{g}_i^{\min} - \varepsilon_i$. Due to $\varepsilon_i > 0$ and (30), the total consumer benefit at $\hat{g}_i$ is below the total production cost of the unit. Hence, the aggregate demand curve at volume $\hat{g}_i$ has prices below $p$ only, which contradicts the assumption $p' \geq p$. In the third case, we have $\overline{g}_i^{\min} - \varepsilon_i < \hat{g}_i$ and all points on the modified supply curve of the unit with the volumes not exceeding $\hat{g}_i$ belong to the unit's supply curve for (31), which has a point with price $p$ and volume $\hat{g}_i$.

If the 3-rd and 4-th bullets of Proposition 4 are applicable, then $\hat{g}_i \leq \overline{g}_i^{\max}$ implies $\hat{g}_i = 0$, and the supply curve of the unit for (31) has zero volumes for any price value.

If conditions of the 5-th bullet of Proposition 4 are satisfied, then the assumption on the aggregate demand curve having a point with the volume $q$ yields $\overline{g}_i^{\min} \leq \hat{g}_i$, which entails $\overline{g}_i^{\min} - \varepsilon_i < \hat{g}_i$ for $\varepsilon_i > 0$. As a result, all points of the modified supply curve of the unit with volumes not exceeding $\hat{g}_i$ belong to the unit's supply curve for (31).

Thus, we conclude that the assumptions stated in this proposition imply that the supply curve of a unit $i$ for (31) has a point with the price $p$ and the volume $\hat{g}_i$.

Hence, the aggregate supply curve for (31) has a point with the price $p$ and the volume $q$. Since the aggregate demand curve for (31) is given by the aggregate demand curve for (5) truncated at volume $q$, condition $p' \geq p$ entails that the aggregate demand and the aggregate supply curves for (31) intersect at the volume $q$. Thus, the intersection point is the optimal point of the dual to (31) and satisfies the power balance constraint. Therefore, $\{\hat{x}, \hat{d}\}$ is the optimal point of (31) and satisfies (14). The proposition is proved.

Let us make an additional simplification: assume that each consumer is allowed to submit a bid composed of the fixed load segment and/or the price-sensitive one-step segment, so that $B_j(d_j) = b_j d_j$, $d_j \in Y_j$, $Y_j = [d_j^{\min}, d_j^{\max}]$, with a constant price $b_j$, $b_j \geq 0$, $\forall j \in J$. If $d_j^{\min} = d_j^{\max}$, then $b_j$ can be assigned an arbitrary non-negative value. (Inclusion in the consumer benefit functions the constant terms that indicate their values at $d_j^{\min}$ will not affect our analysis.) For simplicity, we assume $b_1 > b_2 > b_3 > \ldots > b_m \geq 0$. We also assume zero minimum output limits for all the generating units. Since $\forall j \in J$ the set $cl[\overline{\Omega}_j^{cons}]$ is bounded, it has the maximum element, which we denote as $\overline{d}_j^{\max}$.

*Proposition 6*: If $g_i^{\min} = 0$, $\varepsilon_i > 0$, $\forall i \in I$, and the aggregate modified supply curve has a point with a price $p$ and a volume $q$, $q \geq \sum_{j \in J} d_j^{\min}$, such that the aggregate demand curve has a point with the volume $q$ and a price $p'$, $p' \geq p$, then
$$\sum_{j: j \in J, b_j < p'} d_j^{\min} + \sum_{j: j \in J, b_j \geq p'} \overline{d}_j^{\max} \geq q.$$

*Proof*. Let $j'$ denote the minimum value of $j \in J$ such that $\sum_{j: j \in J, j > j'} d_j^{\min} + \sum_{j: j \in J, j \leq j'} d_j^{\max} \geq q$. Hence, $b_j \geq p'$ for all $j \leq j'$, and if $j' < m$, then $b_j < p'$ for $j > j'$. If $j' = 1$, then the set $\hat{d}_1 = q - \sum_{j \in J \setminus \{1\}} d_j^{\min}$ and $\hat{d}_j = d_j^{\min}$ for the rest of the consumers; if $j' > 1$, then the set $\hat{d}_j = d_j^{\max}$ for all $j \leq j' - 1$, $\hat{d}_{j'} = q - \sum_{j: j \in J, j \leq j'-1} d_j^{\max} - \sum_{j: j \in J, j > j'} d_j^{\min}$, and $\hat{d}_j = d_j^{\min}$ for the rest of the consumers. Clearly,
$$\hat{d}_j \in \arg \max_{d_j \in Y_j} \pi_j^{cons}(p', d_j), \forall j \in J,$$
and $\sum_{j \in J} \hat{d}_j = q$. Application of Proposition 5 entails that $\hat{d}$ satisfies (14) for some $\hat{x}$. Hence, $\hat{d}_j \leq \overline{d}_j^{\max}$, $\forall j \in J$. From $\hat{d}_j = d_j^{\min}$ for $b_j < p'$, we obtain $\sum_{j: j \in J, b_j < p'} d_j^{\min} + \sum_{j: j \in J, b_j \geq p'} \overline{d}_j^{\max} \geq q$. The proposition is proved.

Now we are ready to conclude that to compute the set of market prices (but not the total uplift) in case under consideration, the original consumer feasible sets $Y_j$ can be used in the modified convex hull pricing method, i.e. knowledge of $\overline{\Omega}_j^{cons}$ is not needed to compute prices in the uninode single-period power systems with zero minimum output limits of all units and price-sensitive segment of each consumer demand given by one-step bid.

*Proposition 7*: Let $g_i^{\min}=0$, $\varepsilon_i >0$, $\forall i \in I$, and each consumer bid is composed of the fixed load component with volume $d_j^{\min}$ and/or the price-sensitive one-step segment with a price $b_j$, then the replacement of $\tilde{\pi}_j^{cons}(p,\rho_j)$ by $\pi_j^{cons}(p)$ in (25) doesn't change the set $\tilde{P}^+(\varepsilon,\rho)$ for $\forall \rho_j >0$, $\forall j \in J$.

*Proof*. Let us look at implications of the assumption that the price-sensitive part of the consumer bid consists of the one-step segment only. Since $\bar{d}_j^{\max}$ is the maximum element of $cl[\bar{\Omega}_j^{cons}]$, $\min(\bar{d}_j^{\max}+\rho_j;d_j^{\max})$ is the maximum element of $cl[\tilde{Y}_j(\rho_j)]$. Also, $d_j^{\min}$ is the minimum element of $\tilde{Y}_j(\rho_j)$ and, therefore, of $cl[\tilde{Y}_j(\rho_j)])$ as well. Therefore, $conv(cl[\tilde{Y}_j(\rho_j)])$, which is the convex hull of $cl[\tilde{Y}_j(\rho_j)]$, is given by $[d_j^{\min},\min(\bar{d}_j^{\max}+\rho_j;d_j^{\max})]$. Let us define $B_j(d_j)$ to be $-\infty$ outside $cl[\tilde{Y}_j(\rho_j)]$, then its concave hull on $conv(cl[\tilde{Y}_j(\rho_j)])$ is given by $b_j d_j$. Hence,

$\partial_{\pm}\tilde{\pi}_j^{cons}(p,\rho_j)=-d_j^{\min}$ for $b_j <p$, $\partial_{\pm}\tilde{\pi}_j^{cons}(p,\rho_j)=-\min(\bar{d}_j^{\max}+\rho_j;d_j^{\max})$ for $b_j >p$,

$\partial_{-}\tilde{\pi}_j^{cons}(p,\rho_j)=-\min(\bar{d}_j^{\max}+\rho_j;d_j^{\max})$ for $b_j =p$, $\partial_{+}\tilde{\pi}_j^{cons}(p,\rho_j)=-d_j^{\min}$ for $b_j =p$.

Likewise,

$$\partial_{\pm}\pi_j^{cons}(p)=-d_j^{\min} \text{ for } b_j <p, \quad \partial_{\pm}\pi_j^{cons}(p)=-d_j^{\max} \text{ for } b_j >p,$$

$$\partial_{-}\pi_j^{cons}(p)=-d_j^{\max} \text{ for } b_j =p, \quad \partial_{+}\pi_j^{cons}(p)=-d_j^{\min} \text{ for } b_j =p.$$

Therefore, $\forall p$ we have:

$$\partial_{-}\tilde{\pi}_j^{cons}(p,\rho_j) \geq \partial_{-}\pi_j^{cons}(p), \quad \partial_{+}\tilde{\pi}_j^{cons}(p,\rho_j) \geq \partial_{+}\pi_j^{cons}(p). \quad (32)$$

First, we prove that $\tilde{P}^+(\varepsilon,\rho)$ belongs to a set of solutions to $\{0\} \in \sum_{j\in J}\partial \pi_j^{cons}(p)+\sum_{i\in I}\partial \tilde{\pi}_i^{prod}(p,\varepsilon_i)$. The set $\tilde{P}^+(\varepsilon,\rho)$ is given by the solutions to $\{0\} \in \sum_{j\in J}\partial \tilde{\pi}_j^{cons}(p,\rho_j)+\sum_{i\in I}\partial \tilde{\pi}_i^{prod}(p,\varepsilon_i)$, which is equivalent to

$$\sum_{j\in J}\partial_{-}\tilde{\pi}_j^{cons}(p,\rho_j)+\sum_{i\in I}\partial_{-}\tilde{\pi}_i^{prod}(p,\varepsilon_i) \leq 0 \leq \sum_{j\in J}\partial_{+}\tilde{\pi}_j^{cons}(p,\rho_j)+\sum_{i\in I}\partial_{+}\tilde{\pi}_i^{prod}(p,\varepsilon_i).$$

Eq. (32) yields

$$\sum_{j\in J}\partial_{-}\pi_j^{cons}(p)+\sum_{i\in I}\partial_{-}\tilde{\pi}_i^{prod}(p,\varepsilon_i) \leq \sum_{j\in J}\partial_{-}\tilde{\pi}_j^{cons}(p,\rho_j)+\sum_{i\in I}\partial_{-}\tilde{\pi}_i^{prod}(p,\varepsilon_i).$$

Hence, we need to show that $0 \leq \sum_{j\in J}\partial_{+}\tilde{\pi}_j^{cons}(p,\rho_j)+\sum_{i\in I}\partial_{+}\tilde{\pi}_i^{prod}(p,\varepsilon_i)$ entails $0 \leq \sum_{j\in J}\partial_{+}\pi_j^{cons}(p)+\sum_{i\in I}\partial_{+}\tilde{\pi}_i^{prod}(p,\varepsilon_i)$. These inequalities are equivalent to

$$\sum_{j:j\in J,b_j \leq p} d_j^{\min} + \sum_{j:j\in J,b_j >p}\min(\bar{d}_j^{\max}+\rho_j;d_j^{\max}) \leq \sum_{i\in I}\partial_{+}\tilde{\pi}_i^{prod}(p,\varepsilon_i) \quad (33)$$

and $\sum_{j:j\in J,b_j \leq p}d_j^{\min}+\sum_{j:j\in J,b_j >p}d_j^{\max} \leq \sum_{i\in I}\partial_{+}\tilde{\pi}_i^{prod}(p,\varepsilon_i)$, respectively. Assume the contrary: let the former is true, while the latter is not, then

$$\sum_{j:j\in J,b_j \leq p}d_j^{\min} + \sum_{j:j\in J,b_j >p}d_j^{\max} > \sum_{i\in I}\partial_{+}\tilde{\pi}_i^{prod}(p,\varepsilon_i), \quad (34)$$

which implies $\sum_{j:j\in J,b_j >p}\min(\bar{d}_j^{\max}+\rho_j;d_j^{\max}) < \sum_{j:j\in J,b_j >p}d_j^{\max}$. As a result, for some $\hat{j}\in J$ we have $b_{\hat{j}}>p$ and $\bar{d}_{\hat{j}}^{\max}+\rho_{\hat{j}}<d_{\hat{j}}^{\max}$. Since for any $i$ the point with the price $p$ and the volume $\partial_{+}\tilde{\pi}_i^{prod}(p,\varepsilon_i)$ belongs to the modified supply curve of a generator $i$, the point with the price $p$ and the volume $\sum_{i \in I} \partial_+ \tilde{\pi}_i^{prod}(p, \varepsilon_i)$ belongs to the aggregate modified supply curve. Also, due to $d_j^{min} \leq \min(\bar{d}_j^{max} + \rho_j; d_j^{max})$, $\forall j \in J$, (33) entails $\sum_{j \in J} d_j^{min} \leq \sum_{i \in I} \partial_+ \tilde{\pi}_i^{prod}(p, \varepsilon_i)$. Eq. (34) implies that at the price $p$ the aggregate demand exceeds the aggregate modified supply, which is no lower than the fixed load part of demand. Let us select a point on the aggregate demand curve with volume $\sum_{i \in I} \partial_+ \tilde{\pi}_i^{prod}(p, \varepsilon_i)$, such a point exists because of (34) and continuity of the aggregate demand curve at volumes in the range $[\sum_{j \in J} d_j^{min}, \sum_{j \in J} d_j^{max}]$. If the point is not unique, we take the one with the lowest price. Let $p'$ denote a price at that point, then (34) implies $p' > p$. Proposition 6 gives $\sum_{j: j \in J, b_j < p'} d_j^{min} + \sum_{j: j \in J, b_j \geq p'} \bar{d}_j^{max} \geq \sum_{i \in I} \partial_+ \tilde{\pi}_i^{prod}(p, \varepsilon_i)$. Due to $d_j^{min} \leq \bar{d}_j^{max}$, $\forall j \in J$, replacements of $d_j^{min}$ by $\bar{d}_j^{max}$ for $j$ with $p < b_j < p'$ yields $\sum_{j: j \in J, b_j \leq p} d_j^{min} + \sum_{j: j \in J, b_j > p} \bar{d}_j^{max} \geq \sum_{i \in I} \partial_+ \tilde{\pi}_i^{prod}(p, \varepsilon_i)$, which entails $\sum_{j: j \in J, b_j > p} \min(\bar{d}_j^{max} + \rho_j; d_j^{max}) \leq \sum_{j: j \in J, b_j > p} \bar{d}_j^{max}$. This in turn gives $\sum_{j: j \in J, b_j > p, \bar{d}_j^{max} < d_j^{max}} \min(\bar{d}_j^{max} + \rho_j; d_j^{max}) \leq \sum_{j: j \in J, b_j > p, \bar{d}_j^{max} < d_j^{max}} \bar{d}_j^{max}$ (with the summation set being nonempty subset of $J$ since it includes at least $\hat{j}$), which cannot hold for $\rho_j > 0$, $\forall j \in J$. Thus, we have obtained a contradiction to the assumption. Therefore, $\sum_{j: j \in J, b_j \leq p} d_j^{min} + \sum_{j: j \in J, b_j > p} d_j^{max} \leq \sum_{i \in I} \partial_+ \tilde{\pi}_i^{prod}(p, \varepsilon_i)$, and $\tilde{P}^+(\varepsilon, \rho)$ belongs to a set of solutions to $\{0\} \in \sum_{j \in J} \partial \pi_j^{cons}(p) + \sum_{i \in I} \partial \tilde{\pi}_i^{prod}(p, \varepsilon_i)$.

Now, assume that $\{0\} \in \sum_{j \in J} \partial \pi_j^{cons}(p) + \sum_{i \in I} \partial \tilde{\pi}_i^{prod}(p, \varepsilon_i)$ for some $p$, let us prove that $p \in \tilde{P}^+(\varepsilon, \rho)$. We have

$$\sum_{j \in J} \partial_- \pi_j^{cons}(p) + \sum_{i \in I} \partial_- \tilde{\pi}_i^{prod}(p, \varepsilon_i) \leq 0 \leq \sum_{j \in J} \partial_+ \pi_j^{cons}(p) + \sum_{i \in I} \partial_+ \tilde{\pi}_i^{prod}(p, \varepsilon_i).$$

Using (32), we obtain $0 \leq \sum_{j \in J} \partial_+ \tilde{\pi}_j^{cons}(p, \rho_j) + \sum_{i \in I} \partial_+ \tilde{\pi}_i^{prod}(p, \varepsilon_i)$. Thus, we need to show that $\sum_{j \in J} \partial_- \pi_j^{cons}(p) + \sum_{i \in I} \partial_- \tilde{\pi}_i^{prod}(p, \varepsilon_i) \leq 0$ implies $\sum_{j \in J} \partial_- \tilde{\pi}_j^{cons}(p, \rho_j) + \sum_{i \in I} \partial_- \tilde{\pi}_i^{prod}(p, \varepsilon_i) \leq 0$. Let us assume the contrary, then we have

$$\sum_{i \in I} \partial_- \tilde{\pi}_i^{prod}(p, \varepsilon_i) \leq \sum_{j: j \in J, b_j < p} d_j^{min} + \sum_{j: j \in J, b_j \geq p} d_j^{max}, (35)$$

and $\sum_{i \in I} \partial_- \tilde{\pi}_i^{prod}(p, \varepsilon_i) > \sum_{j: j \in J, b_j < p} d_j^{min} + \sum_{j: j \in J, b_j \geq p} \min(\bar{d}_j^{max} + \rho_j; d_j^{max})$. Since a point with the price $p$ and the volume $\sum_{i \in I} \partial_- \tilde{\pi}_i^{prod}(p, \varepsilon_i)$, which satisfies $\sum_{j \in J} d_j^{min} < \sum_{i \in I} \partial_- \tilde{\pi}_i^{prod}(p, \varepsilon_i) \leq \sum_{j \in J} d_j^{max}$, belongs to the aggregate modified supply curve, let us select a point on the aggregate demand curve with the volume $\sum_{i \in I} \partial_- \tilde{\pi}_i^{prod}(p, \varepsilon_i)$, (if the point is not unique, we take the one with the lowest price), and denote the price at that point as $p'$. Clearly, (35) implies $p' \geq p$. Application of Proposition 6 gives $\sum_{i \in I} \partial_- \tilde{\pi}_i^{prod}(p, \varepsilon_i) \leq \sum_{j: j \in J, b_j < p'} d_j^{min} + \sum_{j: j \in J, b_j \geq p'} \bar{d}_j^{max}$. Therefore, we have $\sum_{i \in I} \partial_- \tilde{\pi}_i^{prod}(p, \varepsilon_i) \leq \sum_{j: j \in J, b_j < p} d_j^{min} + \sum_{j: j \in J, b_j \geq p} \bar{d}_j^{max}$. This, in turn, yields

$\sum_{j:j\in J, b_j \geq p} \min(\bar{d}_j^{\max} + \rho_j; d_j^{\max}) < \sum_{j:j\in J, b_j \geq p} \bar{d}_j^{\max}$, which is impossible for $\rho_j > 0$, $\forall j \in J$. Thus, $\sum_{j\in J} \partial_- \pi_j^{cons}(p) + \sum_{i\in I} \partial_- \tilde{\pi}_i^{prod}(p, \varepsilon_i) \leq 0$ entails $\sum_{j\in J} \partial_- \tilde{\pi}_j^{cons}(p, \rho_j) + \sum_{i\in I} \partial_- \tilde{\pi}_i^{prod}(p, \varepsilon_i) \leq 0$. The proposition is proved.

We emphasize that Proposition 7 implies that if the price-sensitive part of the aggregate demand is formed by the one-step consumer bids only (given the other stated assumptions on the power system), then the prices in the proposed modified convex hull pricing approach can be computed using $Y_j$ - the original feasible sets of the consumers, but to calculate the individual consumer uplifts the sets $\tilde{Y}_j(\rho_j)$ are needed as illustrated in the following example.

*Example* 6.
Consider a power system with one consumer having benefit function $B(d) = bd$ with $b = \$50.00/MWh$, $Y = [0, d^{\max}]$, $d \in Y$, $d^{\max} = 100\,MWh$ and two identical generators with $C_i(x_i) = ag_i + wu_i$, $0 \leq g_i \leq g^{\max}$, $i = 1,2$, $a = \$40.00/MWh$, $w = \$510.00$, $g^{\max} = 80\,MWh$. In this case, the consumer is able to cover the total output cost of just one generator. The primal problem has two solutions: $d^* = 80\,MWh$ with any one of the generators operating at the maximum output with the other generator being offline. Clearly, we have $\bar{\Omega}_i^{prod} = \{(0,0)\} \cup \{(u_i, g_i) | u_i = 1, g_i \in [w/(b-a), g^{\max}]\}$, $i = 1,2$, for the producers and $\bar{\Omega}^{cons} = \{0\} \cup [w/(b-a), g^{\max}]$ for the consumer. For sufficiently small positive $\varepsilon_i$ and $\rho$, we this implies $\tilde{X}_i(\varepsilon_i) = \{(0,0)\} \cup \{(u_i, g_i) | u_i = 1, g_i \in [w/(b-a) - \varepsilon_i, g^{\max}]\}$ and $\tilde{Y}(\rho) = [0, \rho] \cup [w/(b-a) - \rho, g^{\max} + \rho]$. Both the convex hull pricing method and the modified approach produce the same market price $p^+ = \tilde{p}^+(+0,+0,+0) = a + w/g^{\max} = \$46.38/MWh$. Comparison of the pricing outcomes of both methods is summarized in Table 7.

**Table 7.** The uplift payments for Example 6.

|  | Convex hull pricing | | | Modified convex hull pricing | | |
| --- | --- | --- | --- | --- | --- | --- |
|  | $\pi^*$,$ | $\pi^+$,$ | Uplift,$ | $\tilde{\pi}^*$,$ | $\tilde{\pi}^+$,$ | Uplift,$ |
| Consumer | 289.60 | 362.00 | 72.40 | 289.60 | 289.60 | 0.00 |
| Producer 1 | 0.40 | 0.40 | 0.00 | 0.40 | 0.40 | 0.00 |
| Producer 2 | 0.00 | 0.00 | 0.00 | 0.00 | 0.00 | 0.00 |
| Total | 290.00 | 362.40 | 72.40 | 290.00 | 290.00 | 0.00 |

However, if the consumer original feasible set $Y$ is utilized in the modified convex hull pricing method, then we have the same outcome as in the convex hull pricing method. Hence, to deduce the total uplift in the modified convex hull pricing method, the set $\tilde{Y}(\rho)$ cannot be substituted by $Y$.

Now we consider an example illustrating application of the method in case of demand-side non-convexity.

*Example* 7.
Consider a power system with one producer introduced in Example 5 and two consumers having the following benefit functions: $B_1(d_1) = 100d_1$ with $d_1 \in Y_1$, $Y_1 = [0,100]$, and $B_2(y_2) = 80d_2$ with $Y_2 = \{(v_2, d_2) | v_2 = \{0,1\}, d_2 = 200v_2\}$. Clearly, the second consumer has non-convexity associated with opportunity to start an operating cycle that involves consumption of the fixed amount of power, namely, $200 MWh$. Thus, in this example we have both supply-side and demand-side non-convexities. The solution for the primal problem is given by $d_1^* = 50 MWh$, $d_2^* = 200 MWh$, $v_2^* = 1$, $g^* = 250 MWh$, $u^* = 1$. The convex hull pricing method replaces the non-concave benefit function of the consumer 2 with its concave hull $B_2^{**}(y_2) = 80 d_2$, $0 \leq d_2 \leq 200 MWh$, which gives the market price $p^+ = \$80.00/MWh$ and the total uplift of \$1000.00 payable to the consumer 1. Application of the modified approach gives $\overline{\Omega}^{prod} = X$ for the producer, $\overline{\Omega}_1^{cons} = \{0\} \cup \{50\}$ and $\overline{\Omega}_2^{cons} = Y_2$ for the consumers. For any $\varepsilon > 0$, $\rho_2 > 0$ and any sufficiently small $\rho_1 > 0$, we have $\widetilde{X}(\varepsilon) = X$, $\widetilde{Y}_1(\rho_1) = [0, \rho_1] \cup [50 - \rho_1, 50 + \rho_1]$, and $\widetilde{Y}_2(\rho_2) = Y_2$. Application of the modified convex hull pricing results in the same market price $\widetilde{p}^+ = \$80.00/MWh$ but zero total uplift in the limit $\rho_1 \to +0$ since the consumer 1 consumption volumes in the range $(50 + \rho_1, 100]$ are excluded from the lost profit calculation. The outcomes of the both pricing methods are given in Table 8.

**Table 8.** The uplift payments for Example 7.

|  | Convex hull pricing | | | Modified convex hull pricing | | |
|---|---|---|---|---|---|---|
|  | $\pi^*,\$$ | $\pi^+,\$$ | Uplift,$\$$ | $\widetilde{\pi}^*,\$$ | $\widetilde{\pi}^+,\$$ | Uplift,$\$$ |
| Consumer 1 | 1000.00 | 2000.00 | 1000.00 | 1000.00 | 1000.00 | 0.00 |
| Consumer 2 | 0.00 | 0.00 | 0.00 | 0.00 | 0.00 | 0.00 |
| Producer | 14950.00 | 14950.00 | 0.00 | 14950.00 | 14950.00 | 0.00 |
| Total | 15950.00 | 16950.00 | 1000.00 | 15950.00 | 15950.00 | 0.00 |

**6.3 The two-node single-period power system with fixed load**

So far, we have considered uninode single-period power systems only. However, the proposed approach is straightforwardly extended to more complex power systems, which include multiple nodes, power losses, etc., by including the additional constraints in (10) and (14). We give an example that illustrates application of the method for a two-node single-period power system with the network constraint.

*Example* 8.
Consider a system comprised of two nodes connected by a line with transmission capacity of 100 MW. The system has one fixed load consumer and two producers. Both producer 1 and the fixed load of 150 MWh are located in node 1, while producer 2 is connected to a node 2. The market player parameters are given in Table 9.

**Table 9.** Parameter values for Example 8.

|  | $g_i^{min}, MWh$ | $g_i^{max}, MWh$ | $a_i, \$/MWh$ | $w_i, \$$ |
|---|---|---|---|---|
| Producer 1 | 100 | 200 | 15.00 | 20.00 |
| Producer 2 | 150 | 200 | 10.00 | 0.00 |

Clearly, generator 2 is constrained to be offline in the primal problem solution, which is given by $u_1^* = 1$, $g_1^* = 150 MWh$, $u_2^* = g_2^* = 0$. The convex hull pricing method replaces the non-convex cost functions of the producers with the corresponding convex hulls and results in the locational market prices $p_1^+ = \$15.10/MWh$ and $p_2^+ = \$10.00/MWh$ in nodes 1 and 2, respectively. Thus, the congestion price is



$5.10/MWh$. Consequently, the total uplift equals $510.00 and is comprised of payments to producer 1 and FTR holders. The modified approach gives $\overline{\Omega}_1^{prod} = \{(0,150)\}$ for unit 1, and $\overline{\Omega}_2^{prod} = \{(0,0)\} \cup \{(1,0)\}$ for unit 2. For sufficiently small $\varepsilon_1 > 0$, $\varepsilon_2 > 0$, we have $\widetilde{X}_1(\varepsilon_1) = \{(0,0)\} \cup \{(u_1, g_1) \mid u_1 = 1, g_1 \in [150 - \varepsilon_1, 150 + \varepsilon_1]\}$ and $\widetilde{X}_2(\varepsilon_2) = \{(0,0)\} \cup \{(1,0)\}$. Application of the modified convex hull pricing results in the market prices $\widetilde{p}_1^+ = \widetilde{p}_2^+ = \$15.13/MWh$, which imply zero congestion price, and zero total uplift. The outcomes of the both pricing methods are given in Table 10.

**Table 10.** The uplift payments for Example 8.

|  | Convex hull pricing | | | Modified convex hull pricing | | |
|---|---|---|---|---|---|---|
|  | $\pi^*,\$$ | $\pi^+,\$$ | Uplift,$ | $\widetilde{\pi}^*,\$$ | $\widetilde{\pi}^+,\$$ | Uplift,$ |
| Producer 1 | -5.00 | 0.00 | 5.00 | 0.00 | 0.00 | 0.00 |
| Producer 2 | 0.00 | 0.00 | 0.00 | 0.00 | 0.00 | 0.00 |
| FTR holders | 0.00 | 510.00 | 510.00 | 0.00 | 0.00 | 0.00 |
| Total | -5.00 | 510.00 | 515.00 | 0.00 | 0.00 | 0.00 |

The proposed method produces zero congestion rent which manifests the facts that producer 2 is effectively constrained off and there is no power flow between the nodes. Compared to the convex hull pricing, the new approach results in higher market prices, which allows the full compensation of the production cost through the market price without use of the uplift payment.

### 6.4 The uninode two-period power system with price-sensitive load

To illustrate application of the method for a system with time-coupled operation, we study the system with active ramp rate constraint.
*Example* 9.
Consider the uninode two-period power market with one producer and one consumer. Table 11 contains the producer's parameters, which are assumed to be time-independent:

**Table 11.** Parameter values for Example 9.

|  | $g^{min}, MWh$ | $g^{max}, MWh$ | Ramp rate, $MW/hour$ | $a, \$/MWh$ | No-load cost, $\$/hour$ | Start-up cost, $\$$ |
|---|---|---|---|---|---|---|
| Producer | 20 | 100 | 50 | 20.00 | 80.00 | 0.00 |

Initially, the producer's generating unit is assumed to be online operating with output rate of $50MW$. Thus, the producer has the cost function $C(\mathbf{x}) = 20(g^{t=1} + g^{t=2}) + 80(u^{t=1} + u^{t=2})$ and the individual feasible set $X$ given by

$X = \{\mathbf{x} \mid \mathbf{u}g^{min} \leq \mathbf{g} \leq \mathbf{u}g^{max}, \mathbf{u} \in \{0,1\}^2, \mathbf{g} \in R_{\geq 0}^2; -50 \leq g^{t+1} - g^t \leq 50, t = \{0,1\}; g^{t=0} = 50\}$.

The consumer has the fixed load of $80MWh$ in the first time period and both the fixed load of $10MWh$ and price-sensitive load of $30MWh$ in the second time period: $\mathbf{d} = (d^{t=1}, d^{t=2}) \in Y$ with $Y = Y^{t=1} \times Y^{t=2}$, where $Y^{t=1} = \{80\}$ and $Y^{t=2} = [10,40]$. The consumer's price-sensitive component of the benefit function is given by $B(\mathbf{d}) = 10d^{t=2}$. Thus, the ramp-rate constraint provides time coupling between different time periods of a given market planning horizon.

The primal problem solution is given by $\mathbf{u}^* = 1$, $\mathbf{g}^* = \mathbf{d}^* = (80, 30)$. The convex hull pricing method substitutes the non-convex cost function of the producer with the convex hull (which was explicitly constructed in [31] for two-period power system) and gives the market prices $31.60/$MWh$ and $10.00/$MWh$ for $t = 1$ and $t = 2$ respectively with the total uplift of $32.00 payable to the producer.

Application of the proposed pricing method gives the following. Since demand is fixed for $t = 1$, (14) implies that the producer's possible output for $t = 2$ is no lower than $30 MWh$ due to the ramp constraint and is no higher than $30 MWh$ because the marginal benefit of the price-sensitive demand is below the marginal cost of output for $t = 2$ (i.e., $\partial B / \partial d^{t=2} < a$). Consequently, we have $\overline{\Omega}^{prod} = \{(u^{t=1} = 1, g^{t=1} = 80)\} \times \{(u^{t=2} = 1, g^{t=2} = 30)\}$ and $\overline{\Omega}^{cons} = \{d^{t=1} = 80\} \times \{d^{t=2} = 30\}$. Hence, for sufficiently small $\boldsymbol{\varepsilon}$, $\boldsymbol{\rho}$ we have the following modified individual feasible sets of the market players.

$$\widetilde{X}(\boldsymbol{\varepsilon}) = \{(\mathbf{u} = \mathbf{0}, \mathbf{g} = \mathbf{0})\} \cup \{\mathbf{x} \mid \mathbf{x} \in X, 80 - \varepsilon^{t=1} \leq g^{t=1} \leq 80 + \varepsilon^{t=1}, 30 - \varepsilon^{t=2} \leq g^{t=2} \leq 30 + \varepsilon^{t=2}\},$$
$$\widetilde{Y}(\boldsymbol{\rho}) = \{\mathbf{d} = (80, 10)\} \cup \{\mathbf{d} \mid \mathbf{d} \in Y, 30 - \rho^{t=2} \leq d^{t=2} \leq 30 + \rho^{t=2}\}.$$

In the limit of $\boldsymbol{\varepsilon}, \boldsymbol{\rho} \to +0$, the modified pricing scheme yields the market prices $32.67/$MWh$ and $10.00/$MWh$ for $t = 1$ and $t = 2$ respectively with zero total uplift. The pricing outcomes of the both methods are given in Table 12 (with the consumer's profit reflecting the price-sensitive component only).

Table 12. The uplift payments for Example 9.

|  | Convex hull pricing | | | Modified convex hull pricing | | |
| --- | --- | --- | --- | --- | --- | --- |
|  | $\pi^*$,$ | $\pi^+$,$ | Uplift,$ | $\widetilde{\pi}^*$,$ | $\widetilde{\pi}^+$,$ | Uplift,$ |
| Producer | 468.00 | 500.00 | 32.00 | 553.33 | 553.33 | 0.00 |
| Consumer | 0.00 | 0.00 | 0.00 | 0.00 | 0.00 | 0.00 |

The modified convex hull pricing approach states that in this power system the producer has no opportunity to supply any output other than $80 MWh$ for $t = 1$ and any output other than $30 MWh$ for $t = 2$, (likewise for the consumer). To ensure the non-confiscatory pricing for power, the generator state with offline statuses during the whole market planning horizon (which is feasible given the generator's ramp rate and the initial condition for its output) is included in the modified individual feasible set of the producer. For $t = 1$, this implies distribution of the no-load cost not over the maximum output limit of $100 MWh$ (as it is the case in the convex hull pricing method) but over the reduced maximum output limit of $80 MWh$, which results in both higher market price and lower total uplift (which is calculated using the modified individual feasible sets of the market players).

## 7. Discussion and conclusion

The convex hull pricing method treats each output (consumption) volume allowed by the market player individual operational constraints as possible, although some of these volumes may not be technologically and/or economically feasible. We propose applying a market framework that economically discourages any deviations from the output (consumption) schedule set by the centralized dispatch except for those incurred due to participation in the real-time power balancing. In this framework, taking into account either technologically and/or economically infeasible volumes in the lost profit calculation may result in excessive uplift payments to the market players. For each producer (consumer), we propose identifying both technologically and economically feasible output (consumption) volumes as the set of power volumes that can be obtained as solutions to the centralized dispatch problem with new maximum output and consumption limits of the market players not exceeding the corresponding limits used in the original problem. Each of these sets is further extended to include both the points in some small neighborhoods of every output (consumption) volume from the set (provided that these points also satisfy the market player individual operational constraints) and the points ensuring non-confiscatory pricing for power. In the proposed framework, the resulting sets act as the modified individual feasible sets of the corresponding market players in the modified individual decentralized dispatch optimization problems, which enter the dual problem utilized to calculate the set of market prices and the individual uplifts. Thus, the procedure results in extending the constraint set to include new redundant constraints each depending on the output (consumption) volume of the corresponding market player only. This allows formulating the uniform power pricing mechanism without introduction of the new product/services and the associated prices in Lagrangian relaxation procedure. Also, contrary to the convex hull pricing method, the proposed modified approach entails that to explain the resulting uplift payment to a given market player ISO has to pass to it not only the market price but also information that is needed to identify the market player's modified individual feasible set.

The proposed method can be reformulated to fit the widely-used framework that each market player faces a rule stating how its output revenue (consumption cost) is calculated, and, given this rule, the market player lost profit is calculated using the maximum value of its profit on the feasible set specified by its individual (private) constraint set. For definiteness, let us consider a producer and replace all the additional constraints introduced in the proposed method by extending the generator revenue function to include the corresponding penalty functions. For example, for each additional constraint consider the penalty function in the form of the indicator function of this constraint (which is defined to be equal 0 if the constraint holds and $-\infty$, otherwise). From this point of view, the generator revenue function (before the uplift) is the sum of two terms: the linear term (which is the product of the market price and the generator output) and the sum of the indicator functions for all the additional constraints. This revenue function is utilized to calculate the generator lost opportunity cost, which is compensated by the uplift payment, and the market price is then obtained by minimizing the total uplift. In this approach, no new constraints are explicitly included in the individual feasible set, which is used to calculate the generator lost opportunity cost. Also, the indicator functions, which are the nonlinear terms in the pre-uplift revenue function, vanish at both the primal problem and decentralized dispatch optimal points. Therefore, the resulting generator revenue equals the sum of the linear term and the uplift. Thus, the proposed method, which results in linear pricing with the uplift but introduces additional constraints, can be viewed as a specific mixture of the convex hull pricing and the nonlinear pricing without extending the individual feasible set of the market players. The differences between the nonlinear pricing methods and the proposed pricing approach are the following. The nonlinear pricing generally produces market player specific prices and modifies the revenue function so that the primal problem outcome maximizes the market player profit in the decentralized dispatch problem, while the proposed approach results in a universal market price and individual uplifts (which are needed because the profit maximizing volumes in the individual decentralized dispatch problems are generally different from those in the primal problem solution).

We have showed that the proposed modified convex hull pricing method produces the total uplift no higher than that implied by the convex hull pricing algorithm and a generally different set of the market prices. The mathematical reason for potentially lower total uplift payment is the restriction of market player individual feasible sets used in the dual problem that produces both the set of market prices and the lost opportunity costs of market players. Also, if the centralized dispatch problem becomes convex after optimization over the integral variables, then the proposed approach yields the standard marginal prices and zero total uplift payment.

In case of uninode single-period power systems with fixed load and zero minimum output limits, the sets of both technologically and economically feasible output volumes and the corresponding modified individual feasible sets were explicitly constructed in [28]. For such power systems, the proposed method may result in the different sets of market prices and/or total uplift compared to the convex hull pricing approach only in the presence of at least one LNMGU [28]. For fixed load systems with nonzero minimum output limits we show that (contrary to systems with vanishing generator minimum output limits) the modified pricing method yields a set of market prices and a total uplift payment generally different from those produced by the convex hull pricing approach even if no LNMGU is present in the power system.

We have also illustrated that for power systems with the price-sensitive consumer bids the proposed method may produce market pricing (in particular, the total uplift payment) that depends on the structure of demand.

We established the straightforward procedure to find the sets of generator output volumes that are both technologically and economically feasible for uninode single-period power market with price-sensitive consumer bids and zero minimum output limits of all generators. Also, in such power systems we have showed that if the consumer bids are composed of the fixed load and/or one-step price-sensitive bids, then the consumer original feasible sets (and, therefore, the original demand curves) can be used to calculate the set of market prices but not the uplift payments.

Although we restricted our analysis to the uninode power systems, the proposed approach is readily extended to more complex power systems, which include multiple nodes, power losses, etc., as these additional constraints can be incorporated in (10) and (14).

If the modified individual feasible sets can be determined only for some of the market players, the proposed pricing approach with utilization of the original individual feasible sets for the rest of the market players may still lower the total uplift needed to ensure economic stability of the centralized dispatch outcome compared to the convex hull pricing method.

The incentives effects due to market power in both the convex hull pricing method and the proposed approach are complicated due to multi-part bids and non-convex individual feasible sets of the market players, which may spawn the gaming strategies. On one hand, the modified method utilizes smaller individual feasible sets in the lost profit calculation, which may suppress these strategies. On the other hand, since a market player bid and stated parameters of its generating unit may influence the extra constraints added to the individual feasible set of another market participant and the resulting uplifts, such a procedure may stimulate the market power abuse. These issues are to be addressed in further research on the subject.